\documentclass[a4paper,12pt]{article}
\usepackage{amsmath,amsthm,amssymb,latexsym,epsfig,graphicx}
\title{{\bf Doubly connected V-states for the planar Euler equations}}
\author{{Taoufik Hmidi, Francisco de la Hoz, Joan Mateu and Joan Verdera}}
\usepackage{fullpage}
\newtheorem*{teorem}{Theorem A}
\newtheorem*{teorema}{Theorem B}
\newtheorem{teor}{Theorem}

\newtheorem*{teo}{Crandall-Rabinowitz's Theorem}

\newtheorem{co}{Corollary}
\newtheorem{lemma}[co]{Lemma}

\usepackage[colorlinks=true, pdfstartview=FitV, linkcolor=blue, citecolor=blue,
 urlcolor=blue]{hyperref}
\usepackage{color}

\theoremstyle{definition}

\newtheorem{remark}{Remark}
\newtheorem*{gracies}{Acknowledgements}

\newcommand{\T}{\mathbb{T}}
\newcommand{\C}{\mathbb{C}}

\newcommand{\R}{\mathbb{R}}
\newcommand{\Rea}{\operatorname{Re}}
\newcommand{\Ima}{\operatorname{Im}}
\newcommand{\tpi}{\frac{1}{2 \pi i}}

\begin{document}

\date{}

\maketitle

\begin{abstract}
We prove  existence of doubly connected V-states for the planar Euler
equations which are not annuli. The proof proceeds by bifurcation
from annuli at simple ``eigenvalues". The bifurcated $V$-states we
obtain enjoy a $m$-fold symmetry for some $m\ge 3.$  The existence
of doubly connected $V$-states of strict $2$-fold symmetry remains
open.
\end{abstract}

\section{Introduction}

The Euler system in the plane, which governs the motion of a two
dimensional inviscid incompressible fluid, is equivalent, under mild
smoothness assumptions on the velocity field, to the vorticity
equation
\begin{equation}\label{vorticity}
\begin{cases}
\partial_{t}\omega(z,t) +v(z,t)\cdot \nabla\, \omega(z,t)=0, &z\in \C, \; t > 0,\\
v(z,t)= \nabla^{\bot}\triangle^{-1}\omega(z,t), &
 \\ \omega(z,0) =\omega_{0}(z). &
\end{cases} 
\end{equation}
Here  $v(z,t)=v_1(z,t) +i v_2(z,t)$ is the velocity field at the
point $z \in \C$ and time $t$ and the vorticity is given by the
scalar
$$\omega(z,t)=
\partial_1 v_2(z,t) -\partial_2 v_1(z,t), \quad z \in \C, \quad t \ge 0.$$
The known function $\omega_0(z)$ is the initial condition.
  The Biot-Savart law tells us how to recover
velocity from  vorticity. For a fixed time $t$ one has
$$
v(z,t)= \nabla^{\bot}\triangle^{-1}\omega(z,t), \quad z \in \C,
$$
or, in complex notation,
\begin{equation}\label{vorticitycomplex}
v(z,t)=\frac{i}{2\pi} \int_\C
\frac{\omega(\zeta,t)}{\overline{z}-\overline{\zeta}}\,dA
(\zeta),\quad z \in \C,
\end{equation}
with $dA$ being two dimensional Lebesgue measure.
The first equation in \eqref{vorticity} simply means that the
vorticity is constant along particle trajectories. A convenient
reference for these matters is \cite[Chapter 2]{BM}.

Yudovich Theorem asserts that the vorticity equation has a unique
global solution in the weak sense provided the initial vorticity
$\omega_0$ lies in $L^1 \cap L^\infty$. See, for instance,
\cite[Chapter 8]{BM}. A vortex patch is the solution of
\eqref{vorticity} with initial condition the characteristic function
of a bounded domain $D_0.$ Since the  vorticity is transported along
trajectories, we conclude that   $\omega(z,t)$ is the characteristic
function of a domain $D_t.$  In fact,  $D_t=X(D_0,t)$ is the image
of $D_0$ under the flow. Recall that the flow $X$ is the solution of
the ordinary differential equation
\begin{equation}\label{flux}
\frac{d}{dt}X(z,t)= v(X(z,t),t), \quad X(z,0)=z .
\end{equation}
If $D_0$ is the unit disc, then the particle trajectories are
circles centered at the origin and thus $D_t =D_0, \; t\geq0.$ A
remarkable fact discovered by Kirchhoff is that when the initial
condition is the characteristic function of an ellipse centered at
the origin, then the domain $D_t$ is a rotation of $D_0.$ Indeed,
$D_t = e^{i t \Omega}\,D_0$, where the angular velocity $\Omega$ is
determined by the semi-axis $a$ and $b$ of the initial ellipse
through the formula $\Omega = ab/(a+b)^2.$ See, for instance,
\cite[p.304]{BM} or \cite [p.232]{L}. Kirchhoff's result can also be
checked readily using \eqref{Vdef} below.

A rotating vortex patch or V-state is a domain $D_0$ such that if
$\chi_{D_0}$ is the initial condition of the vorticity equation,
then the region of vorticity $1$ rotates with constant angular
velocity around its center of mass, which we assume to be the
origin. In other words, $D_t = e^{i t \Omega}\,D_0$ or,
equivalently, the vorticity at time $t$ is given by
\begin{equation*}\label{Vstate}
\omega(z,t) = \chi_{D_0}(e^{-i t \Omega} z), \quad z \in \C, \quad t
> 0.
\end{equation*}
Here the angular velocity $\Omega$ is a real number associated with
$D_0.$

Deem and Zabusky \cite{DZ} discovered numerically that there exist
simply connected V-states with $m-$fold symmetry for any integer $m
\ge 2.$ A domain $D_0$ is $m$-fold symmetric if $e^{2 \pi i/m} D_0 =
D_0.$ In other words, if it is invariant by the $m-$th dihedral
group, that is, the set of planar isometries leaving invariant a
regular polygon of $m$ sides. A few years later Burbea \cite{B} gave
an analytic proof by bifurcation at simple ``eigenvalues". See also
\cite{HMV}, where the $C^\infty$ boundary regularity of bifurcated
$V-$states close to the disc of bifurcation was proven.
Incidentally, we mention that whether the boundary  of bifurcated
$V-$states is real analytic is an open question.

In this paper we study doubly connected $V$-states. Recall that a
planar domain is doubly connected if its complement in the Riemann
sphere has two connected components. For example, an annulus is
doubly connected. Because of rotation invariance, it is easy to
ascertain that an annulus is a $V$-state. Indeed, if the annulus is
$$
A= \{z \in \C : b <|z|< 1\},
$$
for some inner radius $b$, $0< b <1,$ then the vector field with
vorticity $\chi_A$ is
$$
v(z)= \frac{i}{2}(z-\frac{b^2}{\overline{z}})\chi_A(z)+
\frac{i}{2}\frac{1-b^2}{\overline{z}}\chi_{\C\setminus A}(z).
$$
The trajectories satisfying \eqref{flux} are clearly circles
centered at the origin. Hence vorticity is conserved along
trajectories and $\omega(z,t)= \chi_A(z)$ is a steady solution to
equation \eqref{vorticity}. Therefore $A$ is a $V$-state rotating
with any angular velocity.

 No other explicit doubly connected $V$-state is known. In
\cite{HMV2} one proved that there do not exist doubly connected
Kirchhoff like examples. In other words, the domain between two
ellipses is a $V$-state only if it is an annulus.\\
Our main result reads as follows (a more detailed statement will be
given later in this section).

\begin{teorem}\label{mr}
There exist doubly connected $V$-states which are not annuli.
\end{teorem}

The proof shows that  there exist, like in the simply connected
case, doubly connected $V$-states with $m$-fold symmetry for any
integer $m \ge 3.$  See figure 1, obtained from a numerical
simulation. It is remarkable that our proof breaks down for $m=2.$
In fact, we do not know if there are $V$-states with strict $2$-fold
symmetry, in the sense that they are $2$-fold symmetric but do not
have a $m$-fold symmetry for an even $m$ larger than $2.$ This is
very likely connected to the non existence of doubly connected
Kirchhoff like examples. The difficulty for $m=2$ is that either the
space of ``eigenfunctions" is two dimensional or it is one
dimensional but the transversality condition in
Crandall-Rabinowitz's theorem fails \cite{CR} (see the statement of
this basic result in section 4 below; the transversality condition
is~(d)).

The proof follows the general scheme of \cite{B} and \cite{HMV}. We
first find a system of two equations, each corresponding to a
boundary component of the patch, which describes doubly connected
$V$-states. Each equation is a differentiated form of Burbea's
equation (3.1) in \cite{B} (see also (13) in \cite{HMV}). This
differentiated form was already found useful in \cite [(53)]{HMV} in
proving boundary regularity of $V$-states. Next step is to use
conformal mapping to transport the system into the unit circle $\T$.
We then consider the Banach space of bounded holomorphic functions
on $\{z \in \C : |z|> 1\}$ with derivative satisfying a H\"{o}lder
condition of order $\alpha$ up to the boundary, and  whose extension
to the unit circle has real Fourier coefficients. Here $\alpha$ is
any number satisfying $0 < \alpha < 1.$ We check the hypothesis of
Crandall-Rabinowitz's Theorem for this Banach space and the
transported system. This requires a lengthy but nice technical work.
In particular, we find all possible ``eigenvalues" of the system,
namely, those values of the bifurcation parameter (which is the
angular velocity of rotation) for which the differential of the
mapping giving the system has a non-zero kernel.

This paper is simpler that \cite{HMV} from the technical point of
view. The reason is that the use of the differentiated form of
Burbea's equation for $V$-states smoothes away technical issues. We
also found a much more direct way to deal with complex functions
having real Fourier coefficients, which was unnecessarily involved
in \cite{HMV}. Although throughout the present paper we work in the
doubly connected context, all our proofs apply to the simply
connected case, as the reader will easily realize.

We close this introductory section by stating a more precise form of
Theorem A.
\begin{teorema}\label{mrprima}
Given $0 <b <1,$ let $m$ be a positive integer such that
$$
 1+{b}^{{m}}-\frac{1-{b}^2}{2} {m}<0.
$$
Then there exists a curve of non-annular $m$-fold symmetric doubly
connected $V$-states bifurcating from the annulus $\{z: b < |z| < 1
\}$ at each of the  angular velocities
$$
\Omega_m = \frac{1-b^2}{4}\pm
\frac{1}{2m}\sqrt{\Big(\frac{m(1-b^2)}{2}-1\Big)^2
-b^{2m}}.
$$
\end{teorema}

A remark on the meaning of $\Omega_m$ is in order. As we showed
before an annulus is a $V$-state rotating with any angular velocity.
The angular velocity plays the role of a bifurcation parameter and
$\Omega_m$ is the ``eigenvalue" at which bifurcation takes place. Remark that for each frequency $$ there are two eigenvalues $\Omega_m$ associated with the $\pm$ signs in the previous formula. The reader will find a discussion on the different behavior  of the V-states bifurcating at each of the two values of $\Omega_m$ in Subsection $9.3$.
Another way of understanding $\Omega_m$ is the following. If the
curve of $V$-states is given by a (continuous) mapping
$$
\xi \in (-\epsilon,\epsilon) \mapsto D(\xi),
$$
where $\epsilon$ is a positive number and $D(\xi)$ is a $V$-state
rotating with angular velocity $\Omega(\xi),$ then
 $D(0)$ is the annulus $\{z: b < |z| < 1 \}$ and $\Omega(0)=\Omega_m.$

 Dritschel found in \cite[(4.1), p. 162]{DR} a similar expression for eigenvalues in studying stability of 
 vortices which are a perturbation of an annulus by an eigenfunction associated with a specific mode.
 His $\sigma(m,a)$ is exactly $m \Omega_{m-1}$ with $b$ replaced by $a$.

\section{The equation of a doubly connected $V$-state}
 Let $D$ be a bounded doubly connected domain  with boundary of class $C^1.$ Equivalently, the boundary of $D$ has two connected components
 which are Jordan curves of class $C^1.$ Call $\Gamma_1$ the exterior curve and $\Gamma_2$ the interior one. The goal of this section is
 to deduce an equation which is equivalent to $D$ being a
 $V$-state. Indeed, the equation can be thought of as  a system of two equations, one for each $\Gamma_j, \;
 j=1,2.$

 Consider the vortex patch with initial condition the characteristic function of the domain
 $D.$ At time $t$ the region of vorticity $1$ is a domain $D_t,$
 which we also assume to be of class $C^1.$ The two closed boundary curves
 of $D_t$ are denoted by $\Gamma_{1,\,t}$ and $\Gamma_{2,\,t}.$ We know
 that the  boundary of $D_t$ is advected by the flow \eqref{flux}.
 It is folklore (see,for instance, \cite{B}, \cite{HMV} and \cite{HMV2}) that this condition can be expressed by the equation
\begin{equation}\label{vortex}
\frac{\partial z}{\partial t} (\alpha,t)\cdot \vec{n} =v
(z(\alpha,t),t) \cdot \vec{n},
\end{equation}
where the dot stands for the scalar product of vectors in the plane
and $z(\alpha,t)$ is a proper parametrization of any of the curves
$\Gamma_{j,\,t}, \ j=1,2.$ By this we mean that $z(\alpha,t)$ is
continuously differentiable in $\alpha$ and $t$ and, for fix $t,$ is
a homeomorphism of the interval of parameters $\alpha$ with the
extremes identified onto the closed curve $\Gamma_{j,\,t}.$ The
interpretation of the left hand side of \eqref{vortex} is the normal
component of the motion of the boundary curve and the right  hand
side is the normal component of the motion of a particle on the
curve. Tangential components do not contribute to the motion of the
boundary and are ignored.

The simplest minded argument for  \eqref{vortex} is as follows. Let
$z(\alpha,t)$ and $w(\beta,t)$ two proper parametrizations of one of
the boundary components of $D_t.$ Then there exists a change of
parameters $\alpha(\beta,t)$ such that $w(\beta,t) =
z(\alpha(\beta,t),t)$ for all $\beta$ and $t.$  Thus
\begin{equation*}\label{para}
\frac{\partial w}{\partial t} (\beta,t)=\frac{\partial z}{\partial
\alpha} (\alpha,t) \frac{\partial\alpha (\beta,t)}{\partial t}
+\frac{\partial z}{\partial t} (\alpha,t).
\end{equation*}
Since $\frac{\partial z}{\partial\alpha}(\alpha,t)$ is a tangent
vector at the boundary at the point $z(\alpha,t)$ and
$\frac{\partial \alpha}{\partial t}(\beta,t)$ is a scalar we
conclude that
\begin{equation}\label{boundary}
\frac{\partial w}{\partial t} (\beta,t)\cdot \vec{n} =\frac{\partial
z}{\partial t} (\alpha,t) \cdot \vec{n}
\end{equation}
where $\vec n$ is the exterior unit normal vector at the point
$z(\alpha,t)= w(\beta,t).$ Now apply~\eqref{boundary} with
$w(\beta,t)$ the lagrangian parametrization, that is,
$$
w(\beta,t) = X(w_0(\beta),t),
$$
where $X(z,t)$ is the flow \eqref{flux} and $w_0(\beta)$ is any
proper parametrization of one of the boundary components of $D.$

Let us add to the vortex patch condition \eqref{vortex} the
$V$-state requirement that $D_t$ is a rotation of $D$ around its
center of mass, which we assume to be the origin. This amounts to
say that if $z_0(\alpha)$ is a proper parametrization of one of the
boundary components of $D,$ then $z(\alpha,t)= e^{i \Omega t}
z_0(\alpha)$ is a proper parametrization of the corresponding
boundary component of $D_t.$  Since the scalar product of the
vectors $z$ and $w$ in the plane is just the real part of $z
\overline{w}$, \eqref{vortex} yields
\begin{equation}\label{V1}
\operatorname{Re}(-i \,\Omega \,\overline{z(\alpha,t)}\, \vec{n}) =
\Rea(\overline{v(z(\alpha,t),t)}\, \vec{n}),
\end{equation}
which can be rewritten without resorting to parametrizations as
\begin{equation*}\label{V2}
\operatorname{Re}(-i \,\Omega \,\overline{z}\, \vec{n}) =
\Rea(\overline{v(z,t)}\, \vec{n}), \quad z \in \partial D_t,
\end{equation*}
where $\vec{n}$ is the exterior unit normal vector to the boundary
of $D_t$ at the point $z.$

 By the Biot-Savart law
\eqref{vorticitycomplex}
$$
\overline{v(z,t)} = -\frac{i}{2 \pi} \int_{D_t}
\frac{dA(\zeta)}{z-\zeta}, \quad z \in \partial D_t
$$
and by Green-Stokes
$$
-\frac{1}{\pi} \int_{D_t} \frac{dA(\zeta)}{z-\zeta} = \frac{1}{2 \pi
i} \int_{ \partial D_t} \frac{\overline{\zeta}-
\overline{z}}{\zeta-z}\, d\zeta, \quad z \in D_t.
$$
The last identity remains true also for $z \in \partial D_t,$
because both sides are continuous functions of $z \in \C.$ Therefore
\begin{equation}\label{V3}
\Rea \left[ \left(2\Omega\, \overline{z}\, + \frac{1}{2 \pi i}
\int_{
\partial D_t} \frac{\overline{\zeta}- \overline{z}}{\zeta-z}\,
d\zeta \right)\,\vec{\tau} \right] = 0 , \quad z \in
\partial D_t,
\end{equation}
$\vec{\tau}$ being the unit tangent vector to $\partial D_t,$
positively oriented.

Notice that the left hand side of the above identity is invariant by
rotations. Hence \eqref{V3} holds if and only if it holds for $t=0.$
We conclude that the domain $D$ is a $V$-state if and only if
\begin{equation}\label{Vdef}
\Rea \bigl[ \left(2\Omega\, \overline{z}\, + I(z)\right)\,\vec{\tau}
\bigr] = 0 , \quad z \in
\partial D,
\end{equation}
where
\begin{equation*}\label{I}
I(z)= \frac{1}{2 \pi i} \int_{ \partial D} \frac{\overline{\zeta}-
\overline{z}}{\zeta-z}\, d\zeta, \quad z \in \C.
\end{equation*}

A final remark is that the argument we have discussed gives that
equation \eqref{Vdef} characterizes $V$-states among domains with
$C^1$ boundary, regardless of the number of boundary components. If
the domain is simply connected there is only one boundary component
and so only one equation. If the domain is doubly connected
\eqref{Vdef} gives actually two equations, one per each boundary
component. Of course, in each equation the other boundary component
is present through the operator $I(z).$

\section{Conformal mapping}
In this section we transform \eqref{Vdef} in a system of two
equations on the unit circle $\T$. Living on $\T$ has the advantage
that the system can be posed in a Banach space, so that functional
analysis tools become available.

Recall that our doubly connected bounded domain $D$ has two boundary
components $\Gamma_j, \; j=1,2$, which are Jordan curves of class
$C^1.$ Let $D_j$ be the domain enclosed by the Jordan curve
$\Gamma_j.$ Let $\Delta$ denote the open unit disc $\{z \in \C : |z|
< 1\}.$ The domains $\C \setminus \overline{D_j}$ are simply
connected and thus there are conformal mappings $\Phi_j : \C
\setminus \overline{\Delta} \rightarrow \C \setminus \overline{D_j}$
fixing the point at $\infty.$ We can normalize $\Phi_1$ so that its
expansion at $\infty$ has coefficient $1$ in $z,$ namely,
\begin{equation}\label{fi1}
\Phi_1(z)= z + a_0+ \frac{a_1}{z}+ \dots \equiv z+f(z),
\end{equation}
valid for $z$ outside a large disc. Here $f$ plays the role of an
analytic perturbation of the identity. The expansion of $\Phi_2$ at
$\infty$ is
\begin{equation}\label{fi2}
\Phi_2(z)= b z + b_0+ \frac{b_1}{z}+ \dots \equiv b z+g(z),
\end{equation}
where $0< b <1.$ We can assume the coefficient $b$ to be positive by
making a rotation in $\Delta$.  The inequality $b < 1$ follows from
Schwarz Lemma applied to the mapping $1/(\Phi_2^{-1}\circ
\Phi_1)(1/z), \; |z|< 1.$ As before, $g$ should be viewed as an
analytic perturbation of $b z.$

The domain $D$ can be written as
\begin{equation}\label{domain}
 D= D_1\setminus \overline{D_2}= \left(\C \setminus
\overline{\Phi_1(\C\setminus \overline{\Delta})}\right) \cap
\left(\Phi_2(\C \setminus \overline{\Delta})\right).
\end{equation}
 Notice that if
$f=g=0,$ then $D$ is the annulus $\{z \in \C : b <|z|< 1\}.$

Set
\begin{equation}\label{Ij}
I_j(z)= \frac{1}{2 \pi i} \int_{\Gamma_j} \frac{\overline{\zeta}-
\overline{z}}{\zeta-z}\, d\zeta, \quad z \in \C, \; j=1,2,
\end{equation}
where $\Gamma_j$ is oriented in the counterclockwise direction for
$j=1,2.$  Clearly $\Gamma_j$ can be parametrized by $\Phi_j$ on
$\T.$  Here appears a subtle issue related to the smoothness of
$\Phi_j$ and we pause momentarily to discuss it.

Assume that $A$ is a Jordan domain with $C^1$ boundary, that is, a
bounded simply connected domain whose boundary is a $C^1$ Jordan
curve $\Gamma =\partial A.$ It is well known that the conformal
mapping $\Phi$ of $\C \setminus \overline{\triangle}$ into $\C
\setminus \overline{A}$ extends to a homeomorphism of $\T$ onto
$\Gamma$ and that this homeomorphism is not necessarily continuously
differentiable. This is related to the mapping properties of the
conjugation operator on $\T$, concretely to the fact that it does
not preserve $C^1(\T).$ The Kellogg-Warschawski theorem
\cite[Theorem 3.6, p.49]{P} asserts that if $\Gamma$ is of class
$C^{1+\alpha}$ then $\Phi$ is of class $C^{1+\alpha}(\T)$ (see next
section for a precise definition of this space).

Thus we assume throughout the paper that $D$ is a doubly connected
domain with boundary of class $C^{1+\alpha},$ for some $\alpha$
satisfying $0 < \alpha < 1.$ The two boundary components $\Gamma_j,
\; j=1,2$ are then Jordan curves of class $C^{1+\alpha}.$

Coming back to our previous discussion we conclude that $\Gamma_j$
can be parametrized by $\Phi_j$ on $\T$ and that
$\vec{\tau}(\Phi_j(w))= i w \Phi'_j(w), \; |w|=1.$  Notice that the
preceding equation makes sense at all points $w \in \T$  because
$\Phi_j$ is of class $C^{1+\alpha}(\T) \subset C^1(\T).$  Thus,
taking into account that $\Rea(i z)= -\Ima(z),$ the single equation
\eqref{Vdef} is transformed into the system of two equations on $\T$
\begin{equation}\label{system}
\begin{array}{ll}
\Ima \left[ \left(2\Omega \,\overline{\Phi_1(w})+
I_1(\Phi_1(w))-I_2(\Phi_1(w))\right) w\,\Phi'_1(w)\right]=0, \quad
|w|=1\\*[7pt]
\Ima \left[ \left(2\Omega \,\overline{\Phi_2(w})+
I_1(\Phi_2(w))-I_2(\Phi_2(w))\right) w \, \Phi'_2(w)\right]=0, \quad
|w|=1.
\end{array}
\end{equation}
The functions $(I_j)_{j=1}^2$ introduced in \eqref{Ij} take the form
$$
I_j(z)=\frac{1}{2i \pi}\int_{\mathbb{T}}\frac{\overline{\phi_j(w)}-\overline{z}}{\phi_j(w)-z}\phi_j^\prime(w) dw,\quad z\in \C.
$$
It will be useful in later calculations to replace in the preceding
system the angular velocity $\Omega$ by the parameter
$\lambda=1-2\Omega.$  The left hand sides of the two equations in
\eqref{system} can be thought of as functions of $f$ and $g,$ as
defined in \eqref{fi1} and \eqref{fi2}. Define functions
$F_1(\lambda,f,g)$ and $F_2(\lambda,f,g)$ on $\T$ by
\begin{equation}\label{efaj}
F_j(\lambda,f,g)(w)= \Ima \left[ \left( (1-\lambda)
\,\overline{\Phi_j(w})+ I(\Phi_j(w))\right) w\,\Phi'_j(w)\right],
\quad |w|=1,
\end{equation}
and a function $F(\lambda,f,g)$ by
\begin{equation}\label{efa}
F(\lambda,f,g)(w)= \left(F_1(\lambda,f,g)(w),
F_2(\lambda,f,g)(w)\right), \quad |w|=1.
\end{equation}
Hence the system \eqref{system} is equivalent to the single equation
\begin{equation}\label{efaeq}
F(\lambda,f,g)=0.
\end{equation}

Therefore we have shown that if $D $ is a bounded doubly connected
$V$-state of class $C^{1+\alpha}$, then equation \eqref{efaeq} is
satisfied. Conversely, if $f$ and $g$ are appropriate functions in
$C^{1+\alpha}(\T)$, then $\Phi_1(z)=z+f(z)$ and $\Phi_2(z)=bz+g(z)$
can be extended to conformal mappings of $\C \setminus
\overline{\Delta}$ and the domain $D$ defined by \eqref{domain} is a
$V$-state provided \eqref{efaeq} is satisfied. For example, if $f$
is the boundary values of a function analytic in $\C \setminus
\overline{\Delta}$ with Lipschitz norm
\begin{equation}\label{Lip}
\sup\{\frac{|f(z)-f(w)|}{|z-w|}: |z|>1   |w|>1  \} \equiv \delta < 1
\end{equation} then $\Phi_1(z)=z+f(z)$ is conformal on $\C \setminus
\overline{\Delta},$ because
$$
|\Phi_1(z)-\Phi_1(w)| \ge |z-w|-|f(z)-f(w)| \ge
(1-\delta)|z-w|,\quad |z|>1   |w|>1.
$$
Condition \eqref{Lip} is satisfied provided $f$ belongs to the open
unit ball of $C^{1+\alpha}(\T).$  Thus if $f$ and $g$ are boundary
values of analytic functions on  $\C \setminus \overline{\Delta},$
$f$ belongs to the open unit ball of $C^{1+\alpha}(\T)$ and $g$
belongs to the open ball with center $0$ and radius $b$ in
$C^{1+\alpha}(\T),$ then $\Phi_1(z)=z+f(z)$ and $\Phi_2(z)=bz+g(z)$
are conformal on $\C \setminus \overline{\Delta}$ and the domain $D$
defined by \eqref{domain} is a $V$-state provided \eqref{efaeq} is
satisfied.

In the next section we establish the precise conditions one needs to
require to $f$ and $g$ so that $V$-states are produced via
\eqref{efaeq}.

\section{The Banach spaces for Crandall-Rabinowitz's Theorem}
In this section we discuss the Banach spaces involved in our
application of Crandall-Rabinowitz's Theorem. Its original statement
in \cite[p.325]{CR} is included below for the reader's convenience.
For a linear mapping $L$ we let $N(L)$ and $R(L)$ stand for the
kernel and the range of $L$ respectively. If $Y$ is a vector space
and $R$ is a subspace, then $Y/R$ denotes the quotient space.
\begin{teo}\label{CR}
Let $X$, $Y$ be two Banach spaces, $V$ be a neighborhood of $0$ in $X$ and
$$
F\colon (-1,1) \times V\to Y
$$
have the properties
\begin{enumerate}
\item[(a)] $F(\lambda,0)=0$ for any $|\lambda|<1$.
\item[(b)] The partial derivatives $F_{\lambda}$, $F_{x}$ and $F_{\lambda x}$ exist and are continuous.
\item[(c)] $N(F_{x}(0,0))$ and $Y/R(F_{x}(0,0))$ are one-dimensional.
\item[(d)] $F_{\lambda x}(0,0) x_{0}\notin R(F_{x}(0,0))$, where
$$
N(F_{x}(0,0))=\operatorname{span}\{x_{0}\}.
$$
\end{enumerate}
If  $Z$ is any complement of $N(F_{x}(0,0))$ in $X$, then there is a
neighborhood ~$U$ of $(0,0)$ in $\mathbf{R} \times X$, an
interval~$(-a,a)$, and continuous functions $\varphi\colon (-a,a)\to
\mathbf{R}$, $\psi\colon (-a,a)\to Z$ such that  $\varphi(0)=0$,
$\psi(0)=0$ and
$$
F^{-1} (0)\cap U= \Big\{(\varphi(\xi), \xi x_{0}+\xi \psi(\xi)):| \xi|
<a\big\}\cup /big\{(\lambda,0): (\lambda,0)\in U\Big\}.
$$
\end{teo}

We proceed now to define the spaces $X$ and $Y$ to which the above
theorem will be applied.
Let $E$ be a subset of $\mathbb{C}$ and  $0 < \alpha <1.$ We denote
by $C^\alpha(E)$ the space of continuous functions $f$ such that
$$
\|f\|_{C^\alpha(E)}:=\|f\|_{L^\infty}+\|f\|_\alpha
$$
where $\|f\|_{L^\infty}$ stands for the supremum norm of $f$ on $E$
and
\begin{equation*}
\|f\|_\alpha = \sup_{x\neq y\in E} \frac{|f(x)-f(y)|}
{|x-y|^\alpha}\cdot
\end{equation*}

The space $C^{1+\alpha}(\T)$  is the set of continuously
differentiable functions $f$ on the unit circle $\T$ whose
derivatives satisfy a H\"{o}lder condition of order $\alpha,$ endowed
with the norm
$$
\|f\|_{C^{1+\alpha}(\T)}=\|f\|_{L^\infty}+
\|\frac{df}{dw}\|_{L^\infty}+\Bigl\|\frac{df}{dw}\Bigr\|_{\alpha}.
$$

A word on the operator $d/dw$ is in order. For a smooth function $f$
we set
$$
\frac{df}{dw}=-i e^{-i\theta} \frac{df(e^{i\theta})}{d\theta}.
$$
 It will be more convenient in the sequel, in estimating norms in $C^{1+\alpha}(\T),$ to work with $d/dw$ instead of
$d/d\theta.$  This is legitimate because they differ only by a
smooth factor. Notice that we have the identity
\begin{equation*}\label{diff1}\frac{d \{\overline{f}\}}{dw} =
-\frac{1}{{w}^2} \overline{\frac{df}{dw}}.
\end{equation*}
Let $\C_\infty =\C\cup \{\infty\}$ stand for the Riemann sphere (the
one point compactification of $\C$). Let $C^{1+\alpha}_a(\Delta^c)$
be the space of analytic functions on $\C_\infty \setminus
\overline{\triangle}$ whose derivatives satisfy a H\"{o}lder condition
of order $\alpha$ up to $\T.$ This is also the space
$C^{1+\alpha}_a(\T)$ of functions in $C^{1+\alpha}(\T)$ whose
Fourier coefficients of positive frequency vanish. In other words,
$$
C^{1+\alpha}_a(\Delta^c)=C^{1+\alpha}_a(\T)=\{f \in C^{1+\alpha}(\T)
: f(w)= \sum_{n \ge 0} a_n \overline{w}^n, \; |w|=1\}.
$$
Let $C^{1+\alpha}_{ar}(\T)$ be the subspace of $C^{1+\alpha}_a(\T)$
consisting of those functions in  $C^{1+\alpha}_a(\T)$ with real
Fourier coefficients. This requirement is due to the fact that the
``simple eigenvalues" assumption in condition (c) of
Crandall-Rabinobitz's Theorem could not be proved in our context if
we had worked with the full \emph{complex} Banach space
$C^{1+\alpha}_a(\T)$. At the geometric level this assumption implies
that the $V$-states we will find have the real line as axis of
symmetry.

Define the Banach space $X$ as
\begin{equation}\label{X}
X= C^{1+\alpha}_{ar}(\T) \times C^{1+\alpha}_{ar}(\T).
\end{equation}
 Given $0 < b <1$, let
$V$ stand for $B(0,r_0)\times B(0,r_0),$ where $B(0,r_0)$ is the
open ball of center $0$ and radius $r_0= \frac{1}{2} \min (b, 1-b)$
in $C^{1+\alpha}_{ar}(\T).$ From the above discussion is clear that
if $(f,g) \in V,$ then then $\Phi_1(z)=z+f(z)$ and $\Phi_2(z)= b
z+g(z)$ are conformal on $\C \setminus \overline{\Delta},$ the
Jordan curves $\Gamma_j = \Phi_j(\T)$ are of class $C^{1+\alpha}$
and $\Gamma_2$ is in the domain enclosed by $\Gamma_1.$

Set
$$
H=\{h \in C^\alpha(\T): h(e^{i\theta}) = \sum_{n \ge 1} \beta_n
\sin(n\theta), \, \beta_n \in \R, n \ge 1\},
$$
and define $Y$ as
$$
Y=H \times H.
$$
We now have the basic elements in Crandall-Rabinowitz's Theorem :
the Banach spaces $X$ and $Y,$ the function $F$ (defined in
\eqref{efa}) and its domain $\R \times V.$ We have already mentioned
that $F$ is well defined on $\R  \times V,$ because for $(f,g) \in
V$ $\Phi_1(z)=z+f(z)$ and $\Phi_2(z)= b z+g(z)$ are conformal
mappings on $\C \setminus \overline{\Delta}.$ It is rather easy to
show that $F$ maps $\R  \times V$ into $Y.$  Discussing the details
is the goal of the next section.

\section{$F$ maps into $Y$}
Recall that $F$ was defined in \eqref{efa} as $F=(F_1,F_2)$,  where
\begin{equation}\label{ef1}
F_j(\lambda,f,g)(w)= \Ima \left[ \left( (1-\lambda)
\,\overline{\Phi_j(w})+ I(\Phi_j(w))\right) w \,\Phi'_j(w)\right],
\quad |w|=1.
\end{equation}

To show that $F_j \in C^\alpha(\T)$ we observe that there are three
relevant terms in the right hand side of the above identity :
$\Phi'_j(w),$ which is in
 $C^\alpha(\T)$, $\Phi_j(w),$ which is in $C^{1+\alpha}(\T)$, and $I(\Phi_j(w)),$
 which is in $C^{1+\beta}(\T), \; 0 <\beta <\alpha,$ as was shown in
 \cite[equation(61)]{HMV}. Indeed,  the fact that $I(\Phi_j(w))$ is in $C^\alpha(\T), \; 0< \alpha <1 ,$  follows from
 the following simple lemma (\cite[Lemma 4, p.191]{HMV}), which we state now for
future reference.

 \begin{lemma}
 Let $K(w,\tau)$ be a measurable function on $\T \times \T \setminus \{(w,\tau) \in \T \times \T  : w \neq \tau \}$
 satisfying, for some positive constant $C_0$,
 $$|K(w,\tau)| \le C_0, \quad w \neq \tau, $$
 and that for each $\tau \in \T$ the function $w \rightarrow K(w,\tau)$
 is differentiable for $w \neq \tau$ and
 $$
\left|\frac{\partial}{\partial w} K(w,\tau)\right| \le C_0\,
\frac{1}{|w-\tau|}.
 $$
 Then the integral operator
 \begin{equation}\label{inteop}
T\varphi (w) = \int_{|\tau|=1} K(w,\tau)\,\varphi(\tau)\,d\tau,
 \end{equation}
 satisfies
 $$
\|T\varphi\|_\alpha \le C_\alpha \,C_0 \, \|\varphi\|_\infty, \quad
\varphi \in L^\infty(\T), \quad 0  < \alpha < 1,
 $$
 where $C_\alpha$ depends only on $\alpha.$
 \end{lemma}
The proof of the lemma follows from standard arguments (see, for
example, \cite[p.419] {MOV}).

Proving that the image of $F$ lies in $Y$ is now reduced to
ascertaining that the Fourier series expansion of $F_j(\lambda,f,g)$
is of the form $\sum_{n \ge 1} \beta_n \sin(n\theta), \, \beta_n \in
\R, n \ge 1.$ A function $h$ on $\T$ has a Fourier expansion of that
form if and only if
 $$h(w) = \Ima(\sum_{ n \in  \mathbb{Z}} \beta_n {w}^n), \quad w \in \T, $$ with real
coefficients $\beta_n,  n \in \mathbb{Z}.$ Therefore we have to
prove that
\begin{equation}\label{Gj}
G_j(\lambda,f,g)(w): = \left( (1-\lambda) \,\overline{\Phi_j(w})+
I(\Phi_j(w))\right) w \,\Phi'_j(w), \quad |w|=1,
\end{equation}
 has real Fourier coefficients for $j=1,2.$ Notice that a continuous function
 $H$ defined on the circle $\T$ has real Fourier coefficients if and
 only if
\begin{equation*}\label{RFC}
\overline{H(w)}= H(\overline{w}), \quad |w|=1.
\end{equation*}
 Owing to the definition of the space $X$ \eqref{X} the mappings $\Phi_j, \; j=1,2,$ have real Fourier
 coefficients. Hence all terms appearing in the right-hand side of
 \eqref{Gj} have clearly real Fourier coefficients, except, perhaps,
 $I \circ \Phi_j = I_1 \circ \Phi_j -I_2 \circ \Phi_j , \;j= 1,2.$  Let us deal, for example, with $I_j \circ \Phi_j.$  One simply has to write
\begin{equation*}\label{RFC1}
\begin{split}
\overline{(I_j \circ \Phi_j) (w)} & = - \tpi \int_{|\tau|=1}
\frac{\Phi_j(\tau)-\Phi_j(w)}{\Phi_j(\overline{\tau})-\Phi_j(\overline{w})}\,
\Phi_j'(\overline{\tau})\,d \overline{\tau}\\*[5pt]& = \tpi
\int_{|\zeta|=1}
\frac{\Phi_j(\overline{\zeta})-\Phi_j(w)}{\Phi_j(\zeta)-\Phi_j(\overline{w})}\,
\Phi_j'(\zeta)\,d\zeta\\*[5pt]& = \tpi \int_{|\zeta|=1}
\frac{\overline{\Phi_j(\zeta)-\Phi_j(\overline{w})}}{\Phi_j(\zeta)-\Phi_j(\overline{w})}\,
\Phi_j'(\zeta)\,d\zeta\\*[5pt]& = (I_j \circ \Phi_j)(\overline{w}).
\end{split}
\end{equation*}
The other terms are treated similarly.

\section{Differentiability properties of $F(\lambda,f,g)$}
 Recall that $F(\lambda,f,g)$ is the function that gives
the equation of doubly connected $V$-states \eqref{efaeq}. The goal
of this section is to check the differentiability properties of the
function $F(\lambda,f,g)$  required by Crandall-Rabinowitz's
Theorem.   Notice that $F(\lambda,f,g)$ depends linearly on
$\lambda,$ so that we only have to care about the (joint)
differentiability in $(f,g)$ keeping $\lambda$ fixed.
Differentiability is understood in the Fr\'{e}chet sense. By definition,
$F=(F_1,F_2)$ (see \eqref{efa}). Hence we will work with
$F_j(\lambda, f,g), \; j=1,2,$ as defined in \eqref{efaj}. Examining
the definition of $F_j$  one realizes that the only difficult terms
are
\begin{equation*}\label{31-2}
I(\Phi_j(w))= I_1(\Phi_j(w))-I_2(\Phi_j(w)), \quad |w|=1, \quad
j=1,2,
\end{equation*}
and thus we have to show that the four functions $I_k(\Phi_j(w)), \;
j ,k =1,2,$  are continuously differentiable with respect to the
variable $(f,g)$ in the domain $V$. Recall that
\begin{equation}\label{V}
V= B(0,r_0)\times B(0,r_0), \quad r_0= \frac{1}{2} \min (b, 1-b),
\end{equation}
where $B(0,r_0)$ is the open ball of center $0$ and radius $r_0$ in
the Banach space $C^{1+\alpha}_{ar}(\T).$

 Let $G(f,g)$ be a function of $(f,g)$ defined on $V$ and taking values in $Y.$  We now describe a convenient
 way to prove that $G$ is differentiable on $V,$ that later on will be applied to $I_j \circ \Phi_k  , \; j,k=1,2.$
 One first shows the existence of G\^{a}teaux derivatives in certain particular directions. The G\^{a}teaux derivative of $G$ in the direction $(h,0),\;  h \in
C^{1+\alpha}_{ar}(\T)$ at $(f,g)$ is
\begin{equation*}\label{31-3}
d_f G(f,g)(h): = \lim_{t\rightarrow 0} \frac{G(f+th,g)-G(f,g)}{t},
\end{equation*}
where the limit is required to exist in $Y$ (that is, in
$C^\alpha(\T)$). We will eventually show that $d_f G(f,g)$ is the
standard partial derivative $D_f G(f,g),$ but for now we use the
notation involving the lower case $d.$ Then one checks that $d_f
G(f,g)(h)$ is linear and bounded as a function of $h,$ that is, that
$d_f G(f,g) \in L(C^{1+\alpha}_{ar}(\T),Y).$  The next step is to
prove that $d_f G(f,g)$ is continuous as a mapping of $(f,g) \in V$
into the Banach space $L(C^{1+\alpha}_{ar}(\T),Y).$ In particular,
this shows that, for a fixed $g,$ the mapping $f \rightarrow d_f
G(f,g) \in L(C^{1+\alpha}_{ar}(\T),Y)$ is continuous of $f.$ It is a
well-known elementary fact that then the partial derivative $D_f
G(f,g)$ exists for $(f,g) \in V$ and $D_f G(f,g)=d_f G(f,g).$

One argues similarly for the second variable $g$ and shows that the
limit
\begin{equation*}\label{31-4}
d_g G(f,g)(k): = \lim_{t\rightarrow 0} \frac{G(f,g+tk)-G(f,g)}{t},
\end{equation*}
exists in $Y$ for each $k \in C^{1+\alpha}_{ar}(\T),$ that $d_g
G(f,g) \in L(C^{1+\alpha}_{ar}(\T),Y)$ and that $d_g G(f,g)$ is
continuous as a function of $(f,g) \in V$ into
$L(C^{1+\alpha}_{ar}(\T),Y).$ The conclusion is that the partial
derivative $D_g(f,g)$ exists for $(f,g) \in V$  and $D_g(f,g)=d_g
G(f,g).$

Therefore the partial derivatives $D_f G(f,g)$ and $D_g G(f,g)$
exist for $(f,g) \in V$ and they are continuous functions on $V.$
Thus $G(f,g)$ is continuously differentiable on $V$ (\cite[Chapter
VIII, section 9]{D}).

\subsection{Existence of the G\^{a}teaux derivatives of $F(\lambda,f,g)$}
We first compute the G\^{a}teaux derivative $d_f (I_1 \circ
\Phi_1)(f,g)(h)$ of $I_1 \circ \Phi_1$ at $(f,g)$  in the direction
$(h,0), \; h \in C^{1+\alpha}_{ar}(\T).$ To simplify the writing we
introduce the following notation :
\begin{equation*}\label{deltafi}
\Delta \Phi_1 = \Phi_1(\tau)-\Phi_1(w), \quad\quad\quad \Delta h =
h(\tau)-h(w),
\end{equation*}
and
\begin{equation*}\label{qte}
Q(t,\tau,w)= \frac{\overline{\Delta \Phi_1+t  \Delta h}}{\Delta
\Phi_1+t  \Delta h} \,(\Phi_1'(\tau)+t h'(\tau)),
\end{equation*}
where $t$ is a real number that is close enough to $0$ to ensure
that the denominator does not vanish. We claim that
\begin{equation}\label{gatI}
d_f(I_1\circ\phi_1)(f,g)(h)(w) = \tpi \int_{|\tau|=1}
\frac{\partial}{\partial t} Q(0,\tau,w)\,d\tau,
\end{equation}
or, equivalently, that
\begin{equation*}\label{gatII}
 \int_{|\tau|=1} \left( \frac{Q(t,\tau,w)-Q(0,\tau,w)}{t}- \frac{\partial}{\partial t} Q(0,\tau,w)  \right)
 \,d\tau
\end{equation*}
tends to $0$ in $C^\alpha(\T)$ as $t$ tends to $0.$  A
straightforward computation gives
\begin{equation}\label{gatIII}
\begin{split}
\frac{\partial}{\partial t} Q(0,\tau,w) & = - \frac{(\Delta h)
(\overline{\Delta \Phi_1})}{(\Delta \Phi_1)^2}\,
\Phi_1'(\tau)\\*[5pt]&\quad + \frac{\overline{ \Delta h}}{\Delta
\Phi_1} \,\Phi_1'(\tau)\\*[5pt]&\quad + \frac{\overline{\Delta
\Phi_1}}{\Delta \Phi_1}\,h'(\tau),
\end{split}
\end{equation}
which shows that the right hand-side of \eqref{gatI} is linear as a
function of $h.$  Appealing to Lemma 1 we see that this linear
mapping is bounded from $C^{1+\alpha}_{ar}(\T)$ into $C^\alpha(\T).$
But this fact is a consequence of the proof of \eqref{gatI} we are
going to present. Indeed, \eqref{gatI} follows from Lemma 1 applied
to the kernel
\begin{equation*}\label{ker}
 K_t(\tau,w): = \frac{Q(t,\tau,w)-Q(0,\tau,w)}{t}-
\frac{\partial}{\partial t} Q(0,\tau,w),
\end{equation*}
after checking that the constant of $K_t(\tau,w),$ namely,
$$
C_0(t):= \sup_{\tau \neq w} |K_t(\tau,w)|+ \sup_{\tau \neq w}
|\tau-w|\,|\frac{\partial}{\partial w}K_t(\tau,w)|
$$
tends to $0$ with $t.$

If $\tau \neq w,$ then
\begin{equation*}\label{ker1}
K_t(\tau,w) = \frac{1}{t} \int_{0}^t \left(\frac{\partial}{\partial
u} Q(u,\tau,w)-\frac{\partial}{\partial u} Q(0,\tau,w)\right)\,du
\end{equation*}
and thus
\begin{equation}\label{ker2}
|K_t(\tau,w)| \le \sup_{|u|< |t|} |\frac{\partial^2}{\partial u^2}
Q(u,\tau,w)| \;|t|.
\end{equation}
The derivative of $Q(t,\tau,w)$ with respect to $t$ is given by the sum
\begin{equation*}\label{ker3}
\begin{split}
\frac{\partial}{\partial t} Q(t,\tau,w) &= - \Delta h \,
\frac{\overline{\Delta \Phi_1+t  \Delta h}}{(\Delta \Phi_1+t  \Delta
h)^2}\,(\Phi_1'(\tau)+t h'(\tau))\\*[5pt]
&\quad + \frac{\overline{ \Delta
h}}{\Delta \Phi_1+t  \Delta h} \,(\Phi_1'(\tau)+t h'(\tau))+\frac{\overline{\Delta \Phi_1+t \Delta h}}{\Delta \Phi_1+t  \Delta
h} \,h'(\tau)
\\*[5pt]
\end{split}
\end{equation*}
and the second derivative is described by the sum
\begin{equation}\label{ker4}
\begin{split}
\frac{\partial^2}{\partial t ^2} Q(t,\tau,w) &=  2 (\Delta h)^2 \,
\frac{\overline{\Delta \Phi_1+t  \Delta h}}{(\Delta \Phi_1+t  \Delta
h)^3}\,(\Phi_1'(\tau)+t h'(\tau))\\*[5pt]
&\quad - \frac{\Delta h
\,\overline{\Delta h}}{(\Delta \Phi_1+t  \Delta
h)^2}\,(\Phi_1'(\tau)+t h'(\tau))- \Delta h \,
\frac{\overline{\Delta \Phi_1+t  \Delta h}}{(\Delta \Phi_1+t  \Delta
h)^2}\, h'(\tau)\\*[5pt]
&\quad - \frac{\Delta h \overline{ \Delta h}}{(\Delta
\Phi_1+t \Delta h)^2} \,(\Phi_1'(\tau)+t h'(\tau))+\frac{\overline{ \Delta h}}{\Delta \Phi_1+t  \Delta h}
\,h'(\tau)\\*[5pt]
&\quad-\Delta h\, \frac{\overline{\Delta \Phi_1+t \Delta
h}}{(\Delta \Phi_1+t  \Delta h)^2} \,h'(\tau)+ \frac{\overline{\Delta h}}{\Delta \Phi_1+t  \Delta h} \,h'(\tau).
\end{split}
\end{equation}
Each of the seven terms in \eqref{ker4} can be easily estimated by a
constant  $C(f,h)$ depending only on $f$ and $h.$ Here we are taking
$t$ so small that
$$
|\Delta \Phi_1+t  \Delta h| \ge |\tau-w|-r_0 |\tau-w| -t
\|h\|_{C^{1+\alpha}(\T)}|\tau-w| \ge \frac{1}{3}|\tau-w|.
$$
Therefore, by \eqref{ker2},
$$
|K_t(\tau,w)| \le C(f,h) \;|t|,
$$
which means that the first constant of the kernel tends to $0$ with
$t.$

We now argue similarly to get an estimate for the derivative of
$K_t(\tau,w)$ with respect to $w.$ We have
\begin{equation}\label{ker5}
|\frac{\partial}{\partial w} K_t(\tau,w)| \le \sup_{0 < u< t}
|\frac{\partial^2}{\partial u^2}\frac{\partial}{\partial w}
Q(u,\tau,w)| \;|t|
\end{equation}
and
\begin{equation}\label{ker5}
|\frac{\partial^2}{\partial u^2}\frac{\partial}{\partial w}
Q(u,\tau,w)| \le C(f,h) \frac{1}{|\tau-w|},
\end{equation}
for sufficiently small $t.$  For \eqref{ker5} just differentiate
with respect to $w$ in \eqref{ker4} and notice that that the
absolute value of each term one obtains can be estimated by $ C /
|\tau-w|$. The proof of \eqref{gatI} is now complete.

Since $I_1 \circ \Phi_1$ does not depend on $g,$ one easily sees
that
\begin{equation*}\label{33}
d_g(I_1 \circ \Phi_1)(f,g)=0.
\end{equation*}

The G\^{a}teaux derivatives of the remaining functions $I_1\circ \Phi_2,
I_2\circ \Phi_1$ and $I_2\circ \Phi_2$ are shown to exist as bounded
linear operators from $C^{1+\alpha}_{ar}(\T)$ into $C^\alpha(\T)$ in
the same way. We omit the details.

\subsection{Continuity of $D_f F(\lambda,f,g)$ and $D_g F(\lambda,f,g)$}
We first discuss the continuity of $d_f F(\lambda, f,g)$ with
respect to $(f,g).$ Similar arguments apply for the continuity of
 $d_g F(\lambda,f,g).$ As in the previous subsection we present the
 complete details of just one case. The other cases are dealt with via
 straightforward variations of the case considered.

Take $d_f(I_1\circ\phi_1)(f,g)(h)(w),$ which is the integral in
$\tau,$ divided by $2 \pi i,$ of the three terms in \eqref{gatIII}.
Consider, for example, the integral of the third one
\begin{equation*}\label{cont}
T(f,g)(h)(w): = \tpi \int_{|\tau|=1} \frac{\overline{\Delta \Phi_1
}}{\Delta \Phi_1}\,h'(\tau)\,d\tau,
\end{equation*}
where $\Phi_1(w)=w+f(w)$ and $\Delta \Phi_1 =
\Phi_1(\tau)-\Phi_1(w).$ One has to show continuity of $T(f,g)$ at
the point $(f_0,g_0) \in V$ as a mapping from $V$ into
$L(C^{1+\alpha}_{ar}(\T),Y).$ This case is particularly simple
because $T(f,g)$ does not depend on $g.$ Set $\Phi_{1,0}(w)= w
+f_0(w).$ To estimate $T(f,g)(h)(w)-T(f_0,g_0)(h)(w)$ we just add
and subtract inside the integral the term
$$
\frac{\overline{\Delta \Phi_{1,0}}}{\Delta \Phi_1}\,h'(\tau)
$$
to obtain
\begin{equation*}\label{cont2}
\begin{split}
T(f,g)(h)(w)-T(f_0,g_0)(h)(w) &= \tpi \int_{|\tau|=1}
\frac{\overline{\Delta \Phi_1 - \Delta \Phi_{1,0}}}{\Delta
\Phi_1}\,h'(\tau)\,d\tau \\*[5pt]
&\quad+ \tpi \int_{|\tau|=1} \overline{\Delta
\Phi_{1,0}}\,\frac{\Delta \Phi_{1,0}-\Delta \Phi_1}{\Delta \Phi_1
\,\Delta \Phi_{1,0}}\,h'(\tau)\,d\tau \\*[5pt]
&= A(w)+B(w),
\end{split}
\end{equation*}
where the last identity is a definition of the terms $A(w)$ and
$B(w).$ We estimate  $A$ and $B$ in $C^\alpha(\T)$ by Lemma 1. Think
of the integrands of $A(w)$ and $B(w)$ as kernels $K_A(\tau,w)$ and
$K_B(\tau,w),$ so that $A(w)$ and $B(w)$ are the integrals of the
respective kernels in $\tau$ against the bounded function $1.$ The
straightforward estimate of the absolute value of $K_A$ is
$$
|K_A(\tau,w)| \le \|\frac{1}{\Phi_1'}\|_\infty
\,\|f-f_0\|_{C^{1+\alpha}(\T)}\, \|h'\|_\infty.
$$
For the kernel of $B(w)$  we have
$$
|K_B(\tau,w)| \le \|\frac{1}{\Phi_1'}\|_\infty \,
\|\frac{1}{\Phi_{1,0}'}\|_\infty
\,\|f-f_0\|_{C^{1+\alpha}(\T)}\,\|f_0\|_{C^{1+\alpha}(\T)}\,
\|h'\|_\infty.
$$
Since $\|f_0\|_{C^{1+\alpha}(\T)} \le 1$ and $\|1 /\Phi_1'\|_\infty
\le 2$, because of the definition of $V,$ we get
\begin{equation*}\label{cont3}
|K_A(\tau,w)|+|K_B(\tau,w)| \le 6 \,\|f-f_0\|_{C^{1+\alpha}(\T)}\,
\|h\|_{C^{1+\alpha}(\T)}.
\end{equation*}
Similar estimates yield
\begin{equation*}\label{cont4}
|\partial_w K_A(\tau,w)|+|\partial_w K_B(\tau,w)| \le C\,\frac{
\|f-f_0\|_{C^{1+\alpha}(\T)}\, \|h\|_{C^{1+\alpha}(\T)}} {|\tau-w|},
\end{equation*}
where $C$ is a n absolute constant.

Thus, by Lemma 1,
\begin{equation*}\label{cont5}
\|T(f,g)-T(f_0,g_0)\|_{L(C^{1+\alpha}_{ar}(\T),Y)} \le
C\,\|f-f_0\|_{C^{1+\alpha}(\T)}.
\end{equation*}

\subsection{Second order derivatives}
In this subsection we remark that
\begin{equation}\label{der2}
\frac{\partial}{\partial \lambda}DF(\lambda,f,g)
\end{equation}
exists and is a continuous function of its variables.  This is
straightforward because $F(\lambda,f,g)$ depends linearly on
$\lambda.$ We easily get
\begin{equation}\label{der2f1}
\frac{\partial}{\partial
\lambda}DF_1(\lambda,f,g)(h,k)(w)=-\operatorname{Im}\bigl[w\,
\Phi_1'(w) \, \overline{h(w)}+ w \,\overline{\Phi_1(w)}\,
h'(w)\bigr]
\end{equation}
and
\begin{equation}\label{der2f2}
\frac{\partial}{\partial
\lambda}DF_2(\lambda,f,g)(h,k)(w)=-\operatorname{Im}\bigl[w\,
\Phi_2'(w)
 \, \overline{k(w)}+ w\,\overline{\Phi_2(w)}\, k'(w)\bigr].
\end{equation}
It is then clear that \eqref{der2} is a continuous function of
$(\lambda,f,g)\in \R \times X$ into the space of bounded linear
mappings $L(X \times X, Y).$

\section{Spectral study}

By an eigenvalue we understand a real number $\lambda$ such that the
kernel of $DF(\lambda,0,0)$ is non-trivial. Our plan is to apply
Crandall-Rabinowitz's Theorem to the equation of $V$-states
$F(\lambda,f,g)=0.$  Hence we need to perform a spectral study of
the linearized operator at the annular solution $(\lambda,0,0)$. In
particular we shall identify the "eigenvalues" corresponding to
one-dimensional kernels and determine when the linearized operator
is a Fredholm operator of zero index. Since $F=(F_1,F_2),$ given
$(h,k) \in X,$ we have
\begin{eqnarray*}\label{38}
DF(\lambda,0,0)(h,k)&=& \begin{pmatrix} D_f F_1(\lambda,0,0)h+D_g
F_1(\lambda,0,0)k \\*[4pt] D_f F_2(\lambda,0,0)h+D_g
F_2(\lambda,0,0)k
\end{pmatrix}\\
&:=&\mathcal{L}_\lambda(h,k).
\end{eqnarray*}
Before stating the main result of this section we shall introduce the following set describing the dispersion relation.
\begin{equation}
\mathcal{S}=\Big\{\lambda\in \mathbb{R} : \Delta_n(\lambda,b)=0
\;\;\text{for some non-negative integer}\;\; n \Big\},
\end{equation}
with
$$
\Delta_n(\lambda,b):=\Big(\big(1-\lambda    \big)+b^2 +n
(b^2-\lambda)\Big)\Big(n(1-\lambda)-\lambda \Big)+ b^{2n+2}.
$$
The meaning of $\Delta_n(\lambda,b)$ will become clear in
\eqref{DSn}. The implementation of Crandall-Rabinowitz theorem is
connected to the following theorem which is the cornerstone of the
proof of Theorem B.

\begin{teor}\label{pro1}
The  following assertions hold true.
\begin{enumerate}
\item The kernel of $\mathcal{L}_\lambda: X\to Y$ is non trivial if and only if $\lambda\in \mathcal{S}.$
If in addition $\lambda\neq \frac{1+b^2}{2}$ then the kernel is the
one-dimensional vector space generated by
\begin{equation*}
w\in \T\mapsto \bigl(\,(m(1-\lambda)-\lambda)\,\overline{w}^m, -b^{m}
\,\overline{w}^m \bigr),
\end{equation*}
where $m$ the unique integer such that $\Delta_m(\lambda,b)=0.$
\item If $\lambda= \frac{1+b^2}{2},$ then $\textnormal{dim} \operatorname{Ker}\mathcal{L}_\lambda\in\{1,2\}.$
The kernel has dimension $2$ if and only if there exists $n\geq2$
such that $\Delta_n(\frac{1+b^2}{2},b)=0.$

\item For $\lambda\in \mathcal{S}\backslash\big\{1, b^2,\frac{1+b^2}{2}\big\}$  the range of $\mathcal{L}_\lambda$ is closed and is of codimension  one.
\item  For $\lambda\in \big\{1, b^2\big\}$, the codimension of the range  is infinite.
\item For $\lambda=\frac{1+b^2}{2},$ the codimension of the range is $1$ or $2$. It is $2$ if and only if there exists $n\geq2$ such that $\Delta_n(\frac{1+b^2}{2},b)=0.$
\item The transversality assumption is satisfied if and only if $\lambda\in \mathcal{S}\backslash\big\{\frac{1+b^2}{2}\big\}$.
\end{enumerate}

\end{teor}

\begin{remark}
The transversality assumption is automatically satisfied when
$\lambda\in \mathcal{S}$ and the associated wave number $m$ is zero.
However for $m\geq 1,$ since the function $\lambda\mapsto
\Delta_m(\lambda,b)$ is polynomial of degree $2$, the transversality
condition holds if and only if the discriminant is strictly
positive, that is,
$$
1+b^{m+1}-\frac{1-b^2}{2}(1+m)<0.
$$
\end{remark}
The proof of this theorem will be presented in several steps spread
out in several subsections.  The first step is to have at our
disposal an explicit expression for the functions  $F_1$ and $F_2$
which is suitable for the computations one needs to perform to
describe the linearized operator.

\subsection{More explicit expressions for $F_1$ and  $F_2$}

The non-explicit terms in the definition of $F_j$ in \eqref{efaj}
are $I_j(\Phi_k(w)), \; k =1,2.$ For $I_1(\Phi_1(w)),$ set
$\Phi_1(\tau)=\tau+f(\tau)$ and
\begin{equation*}\label{J1}
J_1(\Phi_1(w))= \tpi \int_{|\tau|=1}
\frac{\overline{f(\tau)-f(w)}}{\Phi_1(\tau)-\Phi_1(w)}\,
\Phi_1'(\tau)\,d \tau,\quad |w|=1.
\end{equation*}
We get, using $\overline{\tau-w}= -\overline{w} (\tau-w)/\tau $ for
$ |\tau|=|w|=1,$
\begin{equation*}\label{39}
\begin{split}
I_1(\phi_1(w))&=\tpi\int_{|\tau|=1}\frac{\overline{\tau-w}}{\Phi_1(\tau)-\Phi_1(w)}\Phi_1^\prime(\tau)d\tau+J_1(\Phi_1(w))\\*[5pt]
&= -\overline{w}
\,\tpi\int_{|\tau|=1}\frac{\tau-w}{\Phi_1(\tau)-\Phi_1(w)}\Phi_1^\prime(\tau)\,\frac{d\tau}{\tau}+J_1(\Phi_1(w))\\*[5pt]
&=-\overline{w}+J_1(\Phi_1(w)).
\end{split}
\end{equation*}
To check that the integral in the second line above is $1$ one
should realize that the expansion at $\infty$ of the integrand is
$1/\tau+ a_2/\tau^2+...$
Similarly
\begin{equation*}\label{40}
I_2(\phi_2(w))=-b\,\overline{w}+J_2(\Phi_2(w))
\end{equation*}
where
\begin{equation*}\label{J2}
J_2(\Phi_2(w))= \tpi \int_{|\tau|=1}
\frac{\overline{g(\tau)-g(w)}}{\Phi_2(\tau)-\Phi_2(w)}\,
\Phi_2'(\tau)\,d \tau,\quad |w|=1.
\end{equation*}

For $I_1(\Phi_2(w))$ one sets
\begin{equation*}\label{41}
\widetilde{I_1}(\Phi_2(w))=\tpi
\int_{|\tau|=1}\frac{\overline{f(\tau)}}{\Phi_1(\tau)-\Phi_1(w)}\Phi_1^\prime(\tau)d\tau.
\end{equation*}
We get
\begin{equation*}\label{42}
\begin{split}
I_1(\phi_2(w))&=\tpi\int_{|\tau|=1}\frac{\overline{\Phi_1(\tau)-\Phi_2(w)}}{\Phi_1(\tau)-\Phi_2(w)}\Phi_1^\prime(\tau)d\tau\\*[5pt]
&=-\overline{\Phi_2(w)}+\tpi\int_{|\tau|=1}\frac{\Phi_1^\prime(\tau)}{\Phi_1(\tau)-\Phi_2(w)}\frac{d\tau}{\tau}+\widetilde{I_1}(\Phi_2(w))\\*[5pt]
&=-\overline{\Phi_2(w)}+\widetilde{I_1}(\Phi_2(w)),
\end{split}
\end{equation*}
where in the last identity we used that the integral over the unit
circle vanishes because the integrand has a double zero at $\infty.$

For $I_2(\Phi_1(w))$, one sets $\Phi_2(w)=bw+g(w)$ and
\begin{equation*}\label{43}
\widetilde{I_2}(\Phi_1(w))=\tpi
\int_{|\tau|=1}\frac{\overline{g(\tau)}}{\Phi_2(\tau)-\Phi_1(w)}\Phi_2^\prime(\tau)d\tau.
\end{equation*}
We get
\begin{equation*}\label{44}
\begin{split}
I_2(\phi_1(w))&=\tpi\int_{|\tau|=1}\frac{\overline{\Phi_2(\tau)-\Phi_1(w)}}{\Phi_2(\tau)-\Phi_1(w)}\Phi_2^\prime(\tau)d\tau\\*[5pt]
&=\frac{b}{2\pi i}\int_{|\tau|=1}\frac{\Phi_2^\prime(\tau)}{\Phi_2(\tau)-\Phi_1(w)}\frac{d\tau}{\tau}+\widetilde{I_2}(\Phi_1(w))\\*[5pt]
&\quad -\overline{\phi_1(w)}\tpi \int_{|\tau|=1}\frac{\Phi_2^\prime(\tau)}{\Phi_2(\tau)-\Phi_1(w)}d\tau, \\*[5pt]
&=\frac{b}{2\pi
i}\int_{|\tau|=1}\frac{\Phi_2^\prime(\tau)}{\Phi_2(\tau)-\Phi_1(w)}\frac{d\tau}{\tau}+\widetilde{I_2}(\Phi_1(w))
\end{split}
\end{equation*}
because the winding number of $\Gamma_2 =\Phi_2(\T)$ with respect to
$\Phi_1(w)$ is $0.$ Take $p$ with $|p| > 1$ such that
$\Phi_1(w)=\Phi_2(p).$ Then, by the residue theorem, the factor of
$b$ in the first term above is
\begin{equation*}\label{45}
\begin{split}
\tpi\int_{|\tau|=1}\frac{{\Phi_2^\prime(\tau)}}{\Phi_2(\tau)-\Phi_2(p)} \frac{d\tau}{\tau}&=-\frac1p\\*[5pt]
&=-\frac{1}{\Phi_2^{-1}(\Phi_1(w))}.
\end{split}
\end{equation*}
By \eqref{efaj} we have
\begin{equation*}\label{efajG}
2i F_j(\lambda,f,g) = G_j(\lambda,f,g)-\overline{G_j(\lambda,f,g)},
\quad j=1,2,
\end{equation*}
where
\begin{equation*}\label{gij}
G_j(\lambda,f,g) (w)=  \left( (1-\lambda) \,\overline{\Phi_j(w})+
I(\Phi_j(w))\right) w\,\Phi'_j(w), \quad j=1,2.
\end{equation*}
Therefore
\begin{equation}\label{46}
\begin{split}
G_1(\lambda,f,g)(w)=\biggl(-\lambda
\overline{w}&+(1-\lambda)\overline{f(w)}+J_1(\Phi_1(w))\\
&+\frac{b}{\Phi_2^{-1}(\Phi_1(w))}-\widetilde{I_2}(\Phi_1(w))\biggr)
  w\bigl(1+f^\prime(w)\bigr)
\end{split}
\end{equation}
and
\begin{equation*}\label{47}
G_2(\lambda,f,g)(w)=\biggl((1-\lambda) b
\overline{w}-\lambda\overline{g(w)}-J_2(\Phi_2(w))+\widetilde{I_1}(\Phi_2(w))\biggr)
w\bigl(b+g^\prime(w)\bigr).
\end{equation*}

\subsection{Computation of $DF(\lambda,0,0)$}

Since $F_j$ is the imaginary part of $G_j$ and we have the explicit
expressions \eqref{46} and \eqref{47} for $G_j,$ our plan is to
compute the derivatives with respect to $f$ and $g$ at the point
$(\lambda,0,0)$ of all terms appearing in \eqref{46} and \eqref{47}.
We first show that
\begin{equation*}\label{DJ1}
D_f J_1(\Phi_1(\cdot))(\lambda,0,0)=0.
\end{equation*}
If $h \in C^{1+\alpha}_{ar}(\T),$ then
\begin{equation*}\label{48}
\begin{split}
D_f J_1(\Phi_1(\cdot))(\lambda,0,0)(h)(w)&=\frac{d}{dt} \Big|_{t=0}
\tpi\int_{|\tau|=1}
t\frac{\overline{h(\tau)-h(w)}}{\tau-w+t(h(\tau)-h(w))}\bigl(1+t\,
h^\prime(\tau)\bigr) d\tau\\*[5pt] &=\tpi\int_{|\tau|=1}
\frac{\overline{h(\tau)-h(w)}}{\tau-w}d\tau=0,
\end{split}
\end{equation*}
where the last identity is due to the fact that the integrand is a
bounded analytic function in the unit disc $\{\tau \in \C : |\tau| <
1\}.$

Since $ J_1(\Phi_1(\cdot))$ does not depend on $g$,
\begin{equation*}\label{DJ11}
D_g J_1(\Phi_1(\cdot))(\lambda,0,0)=0.
\end{equation*}
By similar arguments
\begin{equation}\label{DJ2}
D_f J_2(\Phi_2(\cdot))(\lambda,0,0)= D_g
J_2(\Phi_2(\cdot))(\lambda,0,0)=0.
\end{equation}

Next we show that
\begin{equation}\label{DfQ}
D_f \left(\frac{b}{\Phi_2^{-1}\circ
\Phi_1}\right)(\lambda,0,0)(h)(w)= -b^2 \overline{w}^2 h(w), \quad h
\in  C^{1+\alpha}_{ar}(\T), \quad |w|=1
\end{equation}
and
\begin{equation}\label{DgQ}
D_g \left(\frac{b}{\Phi_2^{-1}\circ
\Phi_1}\right)(\lambda,0,0)(k)(w)= b^2 \overline{w}^2
k(\frac{w}{b}), \quad k \in  C^{1+\alpha}_{ar}(\T), \quad |w|=1.
\end{equation}

For \eqref{DfQ}, take $\Phi_1(w)=w+t h(w)$ and $\Phi_2(w)=b w,$ so
that
$$\Phi_2^{-1}(\Phi_1(w))= \frac{1}{b} (w+th(w))$$
 and thus
\begin{equation*}\label{Dfi2}
\frac{d}{dt}\bigg|_{t=0}\frac{b}{w+th(w)}= -b \overline{w}^2h(w).
\end{equation*}

For \eqref{DgQ}, take $\Phi_1(w)= w $ and $\Phi_2(w)=b w+t k(w).$
Set $\psi(w,t)= \Phi_2^{-1}(w).$ Then
$$
w=\Phi_2(\psi(w,t))= b \psi(w,t)+t k(\psi(w,t)).
$$
Taking derivative with respect to $t$ and evaluating at $0$ yields
$$
0=b \frac{\partial \psi}{\partial t}(w,0)+ k(\psi(w,0))
$$
or
$$
\frac{\partial \psi}{\partial t}(w,0) = - \frac{1}{b}
k(\frac{w}{b}).
$$
Hence
\begin{equation*}\label{Dfi3}
\frac{d}{dt}\bigg|_{t=0} \frac{1}{\left(b w+t k(w)\right)^{-1}}=
-\frac{\frac{\partial \psi}{\partial t}(w,0)}{\psi(w,0)^2} = b
\overline{w}^2 k(\frac{w}{b}).
\end{equation*}

Clearly
\begin{equation}\label{49}
D_f\widetilde{I_2}(\Phi_1(\cdot))(\lambda,0,0)=0
\end{equation}
 because $\widetilde{I_2}(\Phi_1(\cdot))$ vanishes if $g=0.$ We also have
\begin{equation*}\label{50}
D_g\widetilde{I_2}(\Phi_1(\cdot))(\lambda,0,0)=0.
\end{equation*}
To see that, let $k \in C^{1+\alpha}_{ar}(\T).$  Then
\begin{equation*}\label{51}
\begin{split}
\frac{d}{dt}\Big|_{t=0}\tpi\int_{|\tau|=1} \frac{t\,
\overline{k(\tau)}}{b\tau+t k(\tau)-w}\bigl(b+t k^\prime(\tau)\bigr)
d\tau
&=\frac{b}{2\pi i}\int_{|\tau|=1} \frac{\overline{k(\tau)}}{b\tau-w}d\tau\\*[5pt]
&=0.
\end{split}
\end{equation*}
The last identity is due to the fact that the integrand is analytic
in the open unit disc and continuous up to the closed unit disc.

On the one hand,
\begin{equation}\label{52}
D_g\widetilde{I_1}(\Phi_2(\cdot))(\lambda,0,0)=0
\end{equation}
because $\widetilde{I_2}(\Phi_1(\cdot))$ vanishes for $f=0.$ On the
other hand, setting $$h(w)= \sum_{n \ge 0} \alpha_n \frac{1}{w^n},
\quad |w|\ge 1,$$ we get
\begin{equation}\label{53}
\begin{split}
D_f \widetilde{I_1}(\Phi_2(\cdot))(\lambda,0,0)(h)(w)&=\frac{d}{dt}\Big|_{t=0}\tpi\int_{|\tau|=1} \frac{t\, \overline{h(\tau)}}{\tau+t\, h(\tau)-bw}\bigl(1+t\, h^\prime(\tau)\bigr) d\tau\\*[5pt]
&=\tpi\int_{|\tau|=1} \frac{\overline{h(\tau)}}{\tau-b\, w}d\tau\\*[5pt]
&=\sum_{n\geq 0}\alpha_n b^n w^n\\*[5pt]
&=\overline{h({w}/{b})}.
\end{split}
\end{equation}

We are now ready to gather all previous calculations to compute
$DF(\lambda,0,0).$ The expression \eqref{46} of $G_1$, \eqref{DfQ},
\eqref{49} and the product rule for differentiation yield
\begin{equation*}\label{54}
\begin{split}
 D_f G_1(\lambda,0,0)(h)(w)&=\Bigl((1-\lambda) \overline{h(w)}-b^2\overline{w}^2 h(w)\Bigr)w+\bigl(-\lambda \overline{w}+\frac{b^2}{w}\bigr)wh^\prime(w)\\*[3pt]
&=(1-\lambda) w\,\overline{h(w)}-b^2\,\overline{w}
\,h(w)+(b^2-\lambda)\,h^\prime(w)
\end{split}
\end{equation*}
and
\begin{equation*}\label{DF1}
D_f F_1(\lambda,0,0)(h)(w) =\operatorname{Im}
\left[-((1-\lambda)+b^2) \overline{w}\, h(w)+ (b^2-\lambda)
h'(w)\right].
\end{equation*}
Similarly
\begin{equation*}\label{DgG1}
D_g G_1(\lambda,0,0)(h)(w) = b^2 \,\overline{w} \,k(\frac{w}{b})
\end{equation*}
and
\begin{equation*}\label{DgF1}
D_g F_1(\lambda,0,0)(h)(w) =\operatorname{Im} \left[ b^2
\,\overline{w} \,k(\frac{w}{b})\right].
\end{equation*}

By \eqref{DJ2} and \eqref{53} we get
\begin{equation*}\label{DfG2}
D_f G_2(\lambda,0,0)(h)(w) = b \,w\,\overline{h(\frac{w}{b})}
\end{equation*}
and
\begin{equation*}\label{DfF2}
D_f F_2(\lambda,0,0)(h)(w) =- \operatorname{Im} \left[ b
\,\overline{w} \,h(\frac{w}{b})\right].
\end{equation*}
Finally, by \eqref{DJ2} and \eqref{52}
\begin{equation*}\label{DgG2}
D_g G_2(\lambda,0,0)(h)(w) = b \left( -\lambda\, w\,
\overline{k(w)}+ (1-\lambda) k'(w) \right)
\end{equation*}
and
\begin{equation*}\label{DgF2}
D_g F_2(\lambda,0,0)(h)(w) = b  \operatorname{Im} \left[  \lambda \,
\overline{w}\,k(w)+ (1-\lambda) k'(w) \right].
\end{equation*}
Therefore{
\begin{eqnarray}\label{55}
\nonumber DF(\lambda,0,0)(h,k)(w)\!&=&\!
\begin{pmatrix}
\textnormal{Im}\Bigl[-\bigl(1-\lambda\!+\!b^2\bigr)\overline{w}\,h(w)\!+\!(b^2-\lambda)\, h^\prime(w)+b^2\,\overline{w}\,k(\frac{w}{b})\Bigr] \\*[9pt]
\textnormal{Im}\Bigl[-b\,\overline{w}\,h(\frac{w}{b})+b
\bigl(\lambda\overline{w}\,k(w)+(1-\lambda)\,k^\prime(w)\bigr)\Bigr]
\end{pmatrix}\\*[5pt]
&\triangleq&\begin{pmatrix} \mathcal{L}_\lambda^1(h,k)(w) \\*[9pt]
\mathcal{L}_\lambda^2(h,k)(w)\end{pmatrix},
\end{eqnarray}
which gives a convenient expression for the linearized operator
$DF(\lambda,0,0).$ To understand its kernel and range it is useful
to expand the components of \eqref{55} in Fourier series.
 Set
$$h(w) = \sum_{n\ge 0} \alpha_n \overline{w}^n
\quad\quad\quad  \text{and} \quad\quad\quad k(w) = \sum_{n\ge 0}
\beta_n \overline{w}^n.$$ }
Then by straightforward computations we
obtain
%
\begin{equation}\label{FS1}
 \mathcal{L}_\lambda^1(h,k)(e^{i\theta})=\sum_{n\ge 0} \Bigl( \bigl((1-\lambda)+b^2 +n (b^2-\lambda)\bigr) \alpha_n - b^{n+2}\beta_n
 \Bigr)\,\sin((n+1)\theta)
\end{equation}
and
\begin{equation}\label{FS2}
 \mathcal{L}_\lambda^2(h,k)(e^{i\theta})=\sum_{n\ge 0} \Bigl( b^{n+1} \alpha_n + b \bigl(n(1-\lambda)-\lambda \bigr) \beta_n
 \Bigr)\,\sin((n+1)\theta).
\end{equation}
Therefore
\begin{equation}\label{eqw1}
DF(\lambda,0,0)(h,k)(e^{i\theta})=\sum_{n\geq0} M_n\begin{pmatrix}
\alpha_n \\
\beta_n
\end{pmatrix}\sin\big((n+1)\theta\big)
\end{equation}
with
\begin{equation*}
M_n:=\begin{pmatrix}
(1-\lambda)+b^2+n(b^2-\lambda)&-b^{n+2}\\
\\
b^{n+1}& b\bigl(n(1-\lambda)-\lambda\bigr)
\end{pmatrix}.
\end{equation*}
This completes the computation of $DF(\lambda,0,0).$

\subsection{The kernel of $DF(\lambda,0,0)$}
Our next goal  is to derive the dispersion relation which gives the
relationship between the wave number $n$ and the angular velocity
$\Omega=\frac{1-\lambda}{2}$  in order to get a non trivial kernel.
This will be easily follow from  \eqref{FS1} and \eqref{FS2}.
Indeed, the couple of functions  $(h,k)$ is in the kernel of
$DF(\lambda,0,0)$ if and only if all Fourier coefficients in
\eqref{FS1} and \eqref{FS2} vanish, namely,
%
\begin{equation}\label{SKer}
\begin{split}
 \bigl((1-\lambda)+b^2 +n (b^2-\lambda)\bigr) \alpha_n -
b^{n+2}\beta_n &= 0 \\*[5pt]
 b^n \alpha_n + \bigl(n(1-\lambda)-\lambda \bigr) \beta_n &=0
\end{split}
\end{equation}
for $n=0,1,2,\dotsc$ Thus, for each non-negative frequency $n$, we have
a linear homogeneous system of two equations in the unknowns
$\alpha_n$ and $\beta_n.$ The determinant of the system \eqref{SKer}
is
\begin{equation}\label{DSn}
\Delta_n =\Delta_n(\lambda,b)= \Bigl(\bigl(1-\lambda    \bigr)+b^2 +n
(b^2-\lambda)\Bigr)\Bigl(n(1-\lambda)-\lambda \Bigr)+ b^{2n+2}.
\end{equation}

Thus the only way the kernel of $DF(\lambda,0,0)$ can be non-trivial
is that for some frequency $m \ge 0$ one has
$\Delta_m(\lambda,b)=0.$ This non-trivial kernel is one dimensional
if and only if
\begin{equation}\label{1dimker}
\Delta_m(\lambda,b)=0 \quad\quad \text{and}\quad \quad
\Delta_n(\lambda,b) \neq 0, \quad 0\le n \neq m.
\end{equation}
In this case a generator of $\operatorname{Ker} DF(\lambda,0,0)$ is
the pair of functions
\begin{equation}\label{genker}
\bigl(\,(m(1-\lambda)-\lambda)\,\overline{w}^m, -b^{m}
\,\overline{w}^m \bigr), \quad w \in \T.
\end{equation}

We pause to discuss the frequencies $m=0$ and $m=1,$ which turn out
to be specially challenging.

\subsection{Eigenvalues associated with the frequencies $m=0,1$}
For $m=0$ the determinant of the system \eqref{SKer} is
\begin{equation*}\label{DS0}
\Delta_0 = \lambda^2-(1+b^2)\lambda+b^2,
\end{equation*}
which vanishes for $\lambda=1$ and $\lambda=b^2.$

For $\lambda=1$ the determinant $\Delta_n$ is
\begin{equation}\label{deltan}
(1-b^2)\Bigl(n-b^2(1+b^2+\dotsb+b^{2(n-1)})\Bigr) \ge n(1-b^2)^2
\end{equation}
and thus $\Delta_n $ does not vanish for $n \ge 1.$ Hence the kernel
of $DF(\lambda,0,0)$ for $\lambda=1$ is one dimensional and is
generated by $(h,k)=(1,1).$  Therefore $\lambda=1$ is a simple
eigenvalue. We will show in subsection 8.1 below that the
codimension of the range of $DF(1,0,0)$ is infinite, so that
Crandall-Rabinowitz's theorem cannot be applied. It is easily seen
that for $\lambda=1$ or, which is the same, for $\Omega=0,$ equation
\eqref{Vdef} is translation invariant. Thus the translations
$\Phi_{\xi1}(z)= z+\xi$ and $\Phi_{\xi2}(z)=b z+\xi $ give obvious
solutions to \eqref{Vdef} : translated annuli.

For $\lambda=b^2$ the determinant $\Delta_n$ is again the left hand
side of \eqref{deltan} and so it does not vanish for $n \ge 1.$ The
kernel of $DF(b^2,0,0)$ is one dimensional and is generated by
$(b^2,1).$ Thus $\lambda=b^2$ is a simple eigenvalue. As in the
previous case, the codimension of the range of $DF(1,0,0)$ is
infinite, so that Crandall-Rabinowitz's theorem cannot be applied
(see subsection 8.1). We do not know if a curve of solutions to
\eqref{Vdef} emanating from the annulus $\{z:b < |z|<1\}$ can be
found. Equivalently, we do not know if a curve of solutions to
\eqref{efaeq} passing through the solution $\Phi_{1}(z)=z,
\Phi_{2}(z)=b z$ exists. The simple candidate
$$
\Phi_{\xi 1}(z)=z+\xi b^2, \quad |z|\ge 1,
$$
and
$$
\Phi_{\xi 2}(z)=b z+\xi, \quad |z|\ge 1,
$$
fails. Here $\xi$ is a small real number that serves as a parameter
for the curve of candidates. The doubly connected candidate
$V$-state we obtain is the region
$$
A(\xi)= \{z \in \C :  |z-\xi b^2|<1  \quad \text{and}\quad  |z-\xi|>
b\}
$$
between two circles, non-concentric if $\xi \neq 0$. The center of
mass of $A(\xi)$ is the origin, but $A(\xi)$ is not a $V$-state if
$\xi \neq 0$ because the two boundary components are circles and it
was shown in \cite{HMV2} that in this case the inner domain is a
$V$-state only if it is an annulus (which is then centered at the
origin).

We discuss now the eigenvalues associated with the frequency $m=1.$
For $m=1$ one gets
\begin{equation*}\label{Delta1}
\Delta_1= 4 \bigl(\lambda-\frac{(1+b^2)}{2}\bigr)^2
\end{equation*}
and so $\lambda= (1+b^2)/2$ is an eigenvalue. We claim that given $n
\ge 2$ there exists a unique value of $b=b_n$ for which
$\Delta_n((1+b^2)/2,b) = 0.$ Hence, for this particular value of
$b,$ \, $\lambda=(1+b^2)/2$ is a double eigenvalue. To see this, we
first compute the determinant of the system at the frequency $n$ for
$\lambda=(1+b^2)/2$ and we obtain
\begin{equation}\label{Delta2}
\Delta_n(\frac{1+b^2}{2},b)= -\left( \frac{1-b^2}{2} n-
\frac{1+b^2}{2}\right)^2+b^{2n+2},
\end{equation}
which vanishes if and only if
\begin{equation}\label{Delta3}
\frac{1-b^2}{2} n- \frac{1+b^2}{2} = \pm{b^{n+1}}.
\end{equation}
The minus sign above gives the equation
\begin{equation*}\label{E}
E := \frac{1-b^2}{2} n- \frac{1+b^2}{2}+b^{n+1}=0.
\end{equation*}
After some algebra
\begin{equation}\label{E2}
\begin{split}
E & = (1-b) \left(\frac{n-1}{2}(1+b)-b^2
(1+b+\dotsb+b^{n-2})\right)\\*[5pt] &\ge (1-b) (n-1)
(\frac{1+b}{2}-b^2) \ge \frac{(1-b)^2}{2} (n-1)
\end{split}
\end{equation}
and so $E$ is different from zero for $n \ge 2.$ Taking the plus
sign in \eqref{Delta3} we get the equation
\begin{equation}\label{E3}
\varphi(b)= \varphi_n(b) : = (1-b^2)n -(1+b^2)-2b^{n+1}= 0.
\end{equation}
The function $\varphi$ takes the positive value $n-1$ at $0$ and the
negative value $-4$ at $1.$ Hence there is at least one zero between
$0$ and $1.$ This zero is unique because $\varphi$ is strictly
decreasing on $(0,1).$  If $b$ is this zero of $\varphi$,  then
 $\lambda=(1+b^2)/2$ is a double eigenvalue, as it was announced.

If $b$ does not belong to the sequence $\{ b \in (0,1):
\varphi_n(b)=0, \; \text{for some} \; n \ge 2 \}$, then $\lambda=
(1+b^2)/2$ is a simple eigenvalue. However we will show in
subsection 8.2 below that the transversality condition (d) in
Crandall-Rabinowitz's Theorem is not satisfied in this case.

To sum up, for the simple eigenvalues $\lambda=1$ and $\lambda=b^2$
associated with the frequency $m=0$ and for the eigenvalue $\lambda=
(1+b^2)/2$ associated with the frequency $m=1,$ all available
criteria for bifurcation fail. We have not been able to decide
whether or not bifurcation is possible using arguments ``ad hoc".
This seems to be a challenging issue, very likely related for $m=1$
to the fact, proven in \cite{HMV2}, that the region enclosed between
two ellipses which are not circles is not a $V$-state.

\subsection{Eigenvalues associated with frequencies $m \ge 2$}

Fix now $m \ge 2$ and assume that $\Delta_m(\lambda,b)=0$ for some
$\lambda \neq (1+b^2)/2.$  We claim that $\lambda$ is a simple
eigenvalue. Assume, to get a contradiction, that
$\Delta_p(\lambda,b)=0$ for an integer $p
> m$.  The determinant $\Delta_m(\lambda,b)$
is a parabola as a function of $\lambda.$ Indeed we have
\begin{equation*}\label{ema}
\frac{\Delta_m(\lambda,b)}{(m+1)^2}=\lambda^2 -(1+b^2)\lambda +
\frac{(1+b^2) m +b^2 m^2+ b^{2m+2}}{(m+1)^2}.
\end{equation*}
This parabola attains its minimum value at $\lambda= (1+b^2)/2.$ If
$\Delta_m(\lambda,b)= \Delta_p(\lambda,b) =0 , \, m < p$ and
$\lambda \neq (1+b^2)/2,$ then the parabolas corresponding to $m$
and $p$ must be the same. Hence the independent terms should be
equal. The independent term as a function of $m$ is
\begin{equation*}\label{gen}
g(m)= \frac{(1+b^2) m +b^2 m^2+ b^{2m+2}}{(m+1)^2}
\end{equation*}
and its derivative is given by
\begin{equation*}\label{genp}
(m+1)^3 \,g'(m)= -(1-b^2)\left(m-\frac{1+b^2}{1-b^2}\right)+ 2
b^{2m+2} \bigl((m+1)\log(b)-1\bigr).
\end{equation*}
Since $\Delta_m(\lambda,b)=0$ and $\lambda \neq (1+b^2)/2$, we have
$\Delta_m((1+b^2)/2,b) < 0.$ By \eqref{Delta2}
\begin{equation}\label{gen2}
\frac{1-b^2}{2} m- \frac{1+b^2}{2} - b^{m+1} > 0
\end{equation}
or
\begin{equation*}\label{gen3star}
E=\frac{1-b^2}{2} m- \frac{1+b^2}{2} + b^{m+1}< 0.
\end{equation*}
This last possibility is excluded by \eqref{E2} with $n$ replaced by
$m$. Thus one has \eqref{gen2} or, in other words,
\begin{equation}\label{gen3}
m \ge \frac{1+b^2}{1-b^2}+ \frac{2 b^{m+1}}{1-b^2} >
\frac{1+b^2}{1-b^2}.
\end{equation}
But this says that $m $ and $p$ lie in an interval where the
function $g$ is strictly decreasing. Hence $g(m) \neq g(p),$ which
is a contradiction.

\subsection{Codimension of the range of $DF(\lambda,0,0)$}
Assume that for some frequency $m \ge 2$, $\Delta_m(\lambda,b)=0$
and $\Delta_n(\lambda,b)\neq 0, \; 0 \le n \neq m.$  By
\eqref{eqw1}, given $h(w) = \sum_{n\ge 0} \alpha_n \overline{w}^n
\;\, \text{and} \;\, k(w) = \sum_{n\ge 0} \beta_n \overline{w}^n$ in
$C^{1+\alpha}_{ar}(\T)$, we have
\begin{equation}\label{DFhk}
DF(\lambda,0,0)(h,k)= (\varphi,\psi),
\end{equation}
where
\begin{equation}\label{fipsi}
\varphi(e^{i\theta}) = \sum_{n\ge 0}  A_n \, \sin
\bigl((n+1)\theta\bigr), \quad \psi(e^{i\theta}) = \sum_{n\ge 0} B_n\,
\sin \bigl((n+1)\theta \bigr),
\end{equation}
and
\begin{equation}\label{57}
\begin{pmatrix}
A_n\\
B_n
\end{pmatrix}=M_n \begin{pmatrix}
\alpha_n \\
\beta_n
\end{pmatrix},\quad M_n = \begin{pmatrix}
(1-\lambda)+b^2+n(b^2-\lambda)&-b^{n+2}\\
\\
b^{n+1}& b\bigl(n(1-\lambda)-\lambda\bigr)
\end{pmatrix}.
\end{equation}

We know that $\varphi, \psi \in C^{\alpha}(\T)$, $A_n, B_n$ are real
and the vector $(A_m,B_m)$ is in the range of $M_m$ (understood as a
linear mapping from $\R^2$ into itself). Conversely, assume that
$\varphi$ and $\psi$ are functions in $C^\alpha(\T)$ with Fourier
series expansions as in \eqref{fipsi} with real $A_n$ and $B_n.$
Assume, furthermore, that the vector $(A_m,B_m)$ is in the range of
$M_m$. We claim that $(\varphi,\psi)$ is in the range of
$DF(\lambda,0,0),$  which, consequently, has codimension $1$ in $Y.$
To prove the claim take $(\alpha_m,\beta_m)$ satisfying \eqref{57}
(with $n$ replaced by $m$) and $(\alpha_n,\beta_n), 0 \le n \neq m,$
given by
\begin{equation}\label{58}
 \begin{pmatrix}
\alpha_n \\
\beta_n
\end{pmatrix}=M_n^{-1} \begin{pmatrix}
A_n \\
B_n
\end{pmatrix}.
\end{equation}
Define $h(w):= \sum_{n\ge 0} \alpha_n \overline{w}^n \;\, \text{and}
\;\, k(w): = \sum_{n\ge 0} \beta_n \overline{w}^n .$  If we can
prove that the functions $h, k$ belong to $  C^{1+\alpha}(\T),$ then \eqref{DFhk} clearly
holds and we are done. Now,  $h \in C^{1+\alpha}(\T)$ if and only if
\begin{equation}\label{deri}
\sum_{n\ge0} n \,\alpha_n \,\overline{w}^{n+1} \in C^\alpha(\T)
\end{equation}
and $k \in C^{1+\alpha}(\T)$ if and only if
\begin{equation}\label{deri1}
\sum_{n\ge 0}\, n\,\beta_n \,\overline{w}^{n+1}  \in C^\alpha(\T).
\end{equation}
We prove \eqref{deri}. For $ 1 \le n \neq m$, \eqref{58} yields
\begin{equation}\label{enaalfa}
n \,\alpha_n = \,\frac{n}{\Delta_n}\, \Big(
\big[n(1-\lambda)-\lambda\big] A_n + b^{n+1} B_n \Big)
\end{equation}
and
\begin{equation}\label{enabeta}
n \,\beta_n = \,\frac{n}{b \Delta_n}\, \Big(-b^{n+1}\,A_n+
\big[(1-\lambda+b^2+n(b^2-\lambda)\big]\,B_n\Big).
\end{equation}
To illustrate the idea of the proof, take first the term in
\eqref{enaalfa} with fastest growth in the numerator, namely,
\begin{equation*}\label{growthfast}
(1-\lambda)\,\frac{n^2}{\Delta_n} \,A_n.
\end{equation*}
The goal is to prove that
\begin{equation}\label{Fmult}
\theta\mapsto \sum_{n \ge 1} \frac{n^2}{\Delta_n} \,A_n
\,\sin\bigl((n+1)\theta\bigr) \in C^\alpha(\T).
\end{equation}
Set
\begin{equation}\label{deltann}
\Delta_n =A n^2 +B n +C +b^{2n+2},
\end{equation}
where $A, B$ and $C$ depend only on $\lambda$ and $b$. We have $A =
(b^2-\lambda)(1-\lambda),$ so that $A \neq 0$ because we are now
excluding the eigenvalues $\lambda=1$ and $\lambda=b^2$
corresponding to the frequency $m=0.$ Then
\begin{equation*}\label{multip}
\frac{n^2}{\Delta_n} = \frac{1}{A} + \gamma_n,
\end{equation*}
with $|\gamma_n| \le C_0 / n, \; n \ge 1, \;$ for a constant $C_0$
independent of $n.$  Set
\begin{equation*}\label{sigma}
\mu(e^{i\theta}) = \sum_{n \ge 1} \gamma_n \,
\sin\bigl((n+1)\theta\bigr), \quad \theta \in \R.
\end{equation*}
Thus
\begin{equation*}\label{igual}
\sum_{n \ge 1} \frac{n^2}{\Delta_n} \,A_n
\,\sin\bigl((n+1)\theta\bigr) = \frac{1}{A} \varphi(\theta)+(\varphi *
\mu)(\theta), \quad \theta \in \R.
\end{equation*}
By Plancherel's identity  $\mu \in L^2(\T) \subset L^1(\T),$ and so
$\varphi*\mu \in C^\alpha(\T),$ which proves \eqref{Fmult}. The
remaining terms from \eqref{enaalfa} are like $\gamma_n$. This
completes the proof of \eqref{deri}, and \eqref{deri1} is proved
similarly. Notice that the same argument applies to the simple
eigenvalues $\lambda= (1+b^2)/2$ associated with the frequency
$m=1.$

Let us consider the case of the eigenvalues  $\lambda=1$ and
$\lambda=b^2$ corresponding to the frequency $m=0.$ The coefficient
$A$ in \eqref{deltann} vanishes but the coefficient $B$ does not.
For $\lambda=1$ the term that grows faster in \eqref{enabeta} is
\begin{equation*}
 \frac{n^2}{b \Delta_n} (b^2- 1)\,B_n \approx -n \frac{1}{b}\,B_n
\end{equation*}
which means that \eqref{deri1} is in $C^\alpha(\T)$ only if $\psi
\in C^{1+\alpha}(\T).$ Therefore the codimension of the range of
$DF(\lambda,0,0)$ is infinite.

 For $\lambda=b^2$ we consider the
term that grows faster \mbox{in \eqref{enaalfa}.} We obtain
\begin{equation*}\label{qdelta} n \alpha_n  \approx
 n A_n, \quad \text{as} \quad n \rightarrow \infty.
\end{equation*}
Then \eqref{deri} is in $C^\alpha(\T)$ only if $\varphi \in
C^{1+\alpha}(\T)$ and again the codimension of the range of
$DF(\lambda,0,0)$ is infinite.

\subsection{The transversality condition}
Assume that $\lambda$ is a simple eigenvalue and that $v_0$ is a
generator of the kernel of $DF(\lambda,0,0).$  Our goal is to
determine in which cases the assumption $(d)$ in
Crandall-Rabinowitz's theorem is satisfied. This assumption is
\begin{equation*}\label{dcr}
\frac{\partial}{\partial \lambda}DF(\lambda,0,0) (v_0) \notin
R(DF(\lambda,0,0)),
\end{equation*}
where $R(L)$ denotes the range of the mapping $L.$

 By \eqref{der2f1} and \eqref{der2f2} we obtain, setting $f=g=0,$
\begin{equation}\label{der2f1ze}
\frac{\partial}{\partial
\lambda}DF_1(\lambda,0,0)(h,k)(w)=-\operatorname{Im}\Big\{w\,
\overline{h(w)}+ h'(w)\Big\}
\end{equation}
and
\begin{equation}\label{der2f2ze}
\frac{\partial}{\partial \lambda}DF_2(\lambda,f,g)(h,k)(w)=-b
\operatorname{Im}\Big\{w\,
 \, \overline{k(w)}+  k'(w)\Big\},
\end{equation}
for all functions $h, k \in C^{1+\alpha}_{ar}(\T).$ Set $h(w) =
\sum_{n\ge 0} \alpha_n \overline{w}^n \;\, \text{and} \;\, k(w) =
\sum_{n\ge 0} \beta_n \overline{w}^n.$ Then the equations  \eqref{der2f1ze} and
\eqref{der2f2ze} become, if $w=e^{i\theta},$
\begin{equation*}\label{der2f1zef}
\frac{\partial}{\partial \lambda}DF_1(\lambda,0,0)(h,k)(w)=- \sum_{n
\ge 0} (n+1)\, \alpha_n \,\sin\bigl((n+1)\theta\bigr)
\end{equation*}
and
\begin{equation*}\label{der2f2zef}
\frac{\partial}{\partial \lambda}DF_2(\lambda,0,0)(h,k)(w)=-b\,
\sum_{n \ge 0} (n+1)\, \beta_n \,\sin \bigl((n+1)\theta\bigr).
\end{equation*}
We know from \eqref{genker} that a generator of the kernel of
$DF(\lambda,0,0)$ is
\begin{equation*}\label{genker2}
v_0 = \begin{pmatrix}
m(1-\lambda)-\lambda \\
-b^{m}
\end{pmatrix}\,\overline{w}^m.
\end{equation*}
Hence
\begin{equation*}\label{dfv0}
\frac{\partial}{\partial \lambda}DF(\lambda,0,0)(v_0)(w)=
-(m+1)\,\begin{pmatrix}
m(1-\lambda)-\lambda \\
-b^{m+1}
\end{pmatrix}\,\sin
\bigl((m+1)\theta\bigr).
\end{equation*}
Therefore the vector $\frac{\partial}{\partial
\lambda}DF(\lambda,0,0)(v_0)$ is in the range of $DF(\lambda,0,0)$
if and only if the vector $\begin{pmatrix}
m(1-\lambda)-\lambda \\
-b^{m+1}
\end{pmatrix}
\in \R^2$ is a scalar multiple of one column of the matrix $M_m$ in
 \eqref{57}, which is equivalent to
\begin{equation}\label{rangekernel}
\Bigl(m(1-\lambda)-\lambda\Bigr)^2 -b^{2(m+1)}=0.
\end{equation}
Notice that this condition holds for $m=1$ and $\lambda=(1+b^2)/2,$
which tells us that the transversality condition in
Crandall-Rabinowitz's theorem fails for the simple eigenvalues
associated with the frequency $m=1.$

For $m=0$ \eqref{rangekernel} gives $\lambda= \pm b,$ which does not
agree with the possible eigenvalues $\lambda=1$ or $\lambda=b^2$
associated with the frequency $m=0.$ Hence only the case of
frequencies $m \ge 2$ is left. We claim that if $m \ge 2,$ then
\eqref{rangekernel} does not hold.  Combining \eqref{rangekernel}
with $\Delta_m(\lambda,b)=0$ (see \eqref{DSn} for
$\Delta_m(\lambda,b)$) we get, by eliminating $b^{2(m+1)}$,
\begin{equation*}\label{rangek1}
(m+1)\Big(1+b^2-2\lambda\Big)\Big(m(1-\lambda)-\lambda\Big)=0.
\end{equation*}
This gives once again in view of \eqref{rangekernel}
$$
b^{m+1}\Big(1+b^2-2\lambda\Big)=0.
$$
Thus $\lambda=\frac{1+b^2}{2},$ which is not the case because
$\lambda$ is a simple eigenvalue associated with a frequency $m \ge
2.$  Summing up, the transversality condition holds for all simple
eigenvalues except for those of the form $\lambda=(1+b^2)/2$
associated with the frequency $m=1.$
\section{Bifurcation at simple eigenvalues}
In this section we complete the proof that Crandall-Rabinowitz's
theorem can be applied to show that bifurcation is possible at
simple eigenvalues associated with frequencies $m\ge 2$. This, of
course, proves Theorem A. Recall that the differentiability
properties of $F$ have been studied in section 6.  Moreover, Theorem
\ref{pro1} ensures that all the required properties of the
linearized operator are satisfied if and only if $\lambda\in
S\backslash\{\frac{1+b^2}{2}\}$ and the associate wave number $m$ is
bigger than two. These condition can be rewritten as
$$
\Delta_m(\lambda,b)=0, \quad \;\;m\geq2, \quad \;\; \lambda\neq
\frac{1+b^2}{2}.
$$
As $\lambda\mapsto \Delta(m,\lambda)$ is polynomial of degree two
the preceding conditions are equivalent to
$$
\Delta_m( \frac{1+b^2}{2},b)<0.
$$
This inequality has been already discussed in \eqref{gen2} and turns
out to be equivalent to
$$
\frac{1-b^2}{2} m-\frac{1+b^2}{2}-b^{m+1}>0,
$$
which in turn is equivalent to
$$
1+b^{m+1}-\frac{1-b^2}{2}(m+1) < 0.
$$
 At this stage we conclude that  for $m\geq2$ and for each simple
solution of $\Delta_m(\lambda,b)=0$ Crandall-Rabinowitz's theorem
can be applied  and therefore we get a bifurcating curve at the
annulus at the values of $\lambda$
\begin{equation*}\label{angular}
\lambda_m^{\pm} = \frac{1+b^2}{2}\pm
\frac{1}{(m+1)}\sqrt{\bigl(\frac{(m+1)(1-b^2)}{2}-1\bigr)^2
-b^{2(m+1)}},
\end{equation*}
which yield the angular velocities
$$
\Omega_m = \frac{1-b^2}{4}\pm \frac{1}{2
(m+1)}\sqrt{\Big(\frac{(m+1)(1-b^2)}{2}-1\Big)^2 -b^{2(m+1)}}.
$$
These two angular velocities correspond to two curves of $V-$states
which bifurcate at the annulus with the same wave number $m$. Each
point different from the annulus in any of these curves is a non
annular doubly connected $V-$state. The goal of next subsection is
to show that these $V-$states enjoy a (m+1)-fold symmetry. Hence the
proof of Theorem B will be completed  replacing $(m+1)$ by $m$.

\subsection{$(m+1)$-fold symmetry of bifurcated $V$-states}
We have proved that we can bifurcate at simple eigenvalues
associated with a frequency $m \ge 2.$  The purpose of this
subsection is to show that the bifurcated $V$-states enjoy a
$(m+1)$-fold symmetry. This is rather simple to prove by changing
the spaces $X$ and $Y$ appropriately. We replace $X$ by
\begin{equation*}\label{Xm}
X_m= C^{1+\alpha}_{arm}(\T) \times C^{1+\alpha}_{arm}(\T),
\end{equation*}
where $C^{1+\alpha}_{arm}(\T)$ is the space of functions $f \in
C^{1+\alpha}_{ar}(\T)$ with Fourier series expansion
\begin{equation*}\label{Xmseries}
\begin{split}
f(w) & = a_m \overline{w}^{m} + a_{2m+1} \overline{w}^{2m+1}
+\dotsb+a_{n(m+1)-1} \overline{w}^{n(m+1)-1} +\dotsb \\ & = w \Bigl(a_m
\overline{w}^{m+1}+ a_{2m+1} \overline{w}^{2(m+1)}+\dotsb+a_{n(m+1)}
\overline{w}^{n(m+1)} +\dotsb\Bigr), \quad w \in \T.
\end{split}
\end{equation*}
As we did before with $X,$ we let $V$ stand for $B(0,r_0)\times
B(0,r_0),$ where $B(0,r_0)$ is the open ball of center $0$ and
radius $r_0= \frac{1}{2} \min (b, 1-b)$ in $C^{1+\alpha}_{arm}(\T).$
If $(f,g) \in V,$ then $\Phi_1(z)=z+f(z)$ and $\Phi_2(z)= b
z+g(z)$ are conformal mappings with a $(m+1)$-fold symmetry. In
fact, for $\Phi_1$ we have,
\begin{equation*}\label{fi1m}
\Phi_1(z)= z\Bigl(1 + \frac{a_{m}}{z^{m+1}}+
\frac{a_{2m+1}}{z^{2(m+1)}}+\dotsb+
\frac{a_{n(m+1)-1}}{z^{n(m+1)}}+\dotsb\Bigr),
\end{equation*}
which yields
\begin{equation}\label{fi1ms}
\Phi_1( e^{i \frac{2 \pi}{m+1}}z)= e^{i \frac{2 \pi}{m+1}}
\Phi_1(z).
\end{equation}
Similarly, for $\Phi_2$ we get
\begin{equation*}\label{fi2m}
\Phi_1(z)= z\Bigl(b + \frac{b_{m}}{z^{m+1}}+
\frac{b_{2m+1}}{z^{2(m+1)}}+\dotsb+
\frac{b_{n(m+1)-1}}{z^{n(m+1)}}+\dotsb\Bigr),
\end{equation*}
and
\begin{equation}\label{fi2ms}
\Phi_2( e^{i \frac{2 \pi}{m+1}}z)= e^{i \frac{2 \pi}{m+1}}
\Phi_2(z).
\end{equation}
Set
$$
H_m=\Big\{h \in C^\alpha(\T); h(e^{i\theta}) = \sum_{n \ge 1} \beta_n
\sin\big(n(m+1)\theta\big), \, \beta_n \in \R, n \ge 1\Big\},
$$
and define $Y_m$ as
$$
Y_m=H_m \times H_m.
$$
We need to check that $F$ as defined in \eqref{efa} maps $X_m$ into
$Y_m.$ For this it is sufficient to ascertain that, given $(f,g) \in
X_m,$ the Fourier series expansion of $F_j(\lambda,f,g),$ as defined
in \eqref{ef1}, is of the form $\displaystyle{\sum_{n \ge 1} \beta_n
\sin(n(m+1)\theta)}, \, \beta_n \in \R, n \ge 1.$ A function $h$ on
$\T$ has a Fourier expansion of the form above if and only if
 $$h(w) = \Ima \bigl(\sum_{ n \in  \mathbb{Z}} \beta_n {w}^{n(m+1)}\bigr), \quad w \in \T, $$ with real
coefficients $\beta_n,  n \in \mathbb{Z}.$ Therefore we have to
prove that for $j=1,2$
\begin{equation*}\label{Gj1}
G_j(\lambda,f,g)(w): = \left( (1-\lambda) \,\overline{\Phi_j(w})+
I(\Phi_j(w))\right) w \,\Phi'_j(w), \quad |w|=1,
\end{equation*}
 has a Fourier series expansion of the type
\begin{equation*}\label{exp}
\sum_{n \in \mathbb{Z}} \beta_n  w^{n(m+1)}, \quad \beta_n \in \R,
\quad n\in \mathbb{Z}.
\end{equation*}
A function $k(w)$ has a Fourier series expansion as above if and
only if
\begin{equation*}\label{fexp}
k(e^{i \frac{2 \pi}{m+1}} w) = k(w), \quad w \in \T.
\end{equation*}
This follows readily for the term $ \overline{\Phi_j(w)} w
\,\Phi'_j(w).$ For the second term $I(\Phi_j(w)) w \,\Phi'_j(w),$
one has to show that
\begin{equation*}\label{simI}
I(\Phi_j( e^{i \frac{2 \pi}{m+1}}w))= e^{-i \frac{2 \pi}{m+1}}
I(\Phi_j(w)),
\end{equation*}
which is easy, just by looking at the integral defining
$I(\Phi_j(w))$ and making a simple change of variables. This
completes the proof that $F$ maps $X_m$ into $Y_m.$

The rest is straightforward. The kernel of $DF(\lambda,0,0)$ is
generated by \eqref{genker}, which is in $X_m.$ Since we are
assuming that $m\ge 2,$ the codimension of the range of
$DF(\lambda,0,0)$ is still $1$ in $Y_m.$ Finally the transversality
condition holds. Therefore we can apply Crandall-Rabinowitz's
Theorem in $X_m$ and $Y_m$ and we get a curve of solutions to
\eqref{efaeq} of the form
$$\xi\in(-\epsilon, \epsilon)\mapsto (\lambda_\xi, f_\xi,g_\xi) \in \R \times
X_m.$$ The conformal mappings provided by $f_\xi$ and $g_\xi$ are of
the form
\begin{equation*}\label{conf1}
\Phi_{1\xi}(z)= z\Bigl(1 + \xi \frac{a_{1}(\xi)}{z^{m+1}}+ \xi
\frac{a_{2}(\xi)}{z^{2(m+1)}}+\dotsb+ \xi
\frac{a_{n}(\xi)}{z^{n(m+1)}}+\dotsb\Bigr)
\end{equation*}
and
\begin{equation*}\label{conf2}
\Phi_{2\xi}(z)= z\Bigl(b + \xi \frac{b_{1}(\xi)}{z^{m+1}}+ \xi
\frac{b_{2}(\xi)}{z^{2(m+1)}}+\dotsb+ \xi
\frac{b_{n}(\xi)}{z^{n(m+1)}}+\dotsb\Bigr).
\end{equation*}
Thus the $V$-state we obtain is, according to \eqref{domain},
\begin{equation*}\label{domain2}
 D_\xi= D_{1\xi}\setminus \overline{D_{2\xi}}= \left(\C \setminus
\overline{\Phi_{1\xi}(\C\setminus \overline{\Delta})}\right) \cap
\left(\Phi_{2\xi}(\C \setminus \overline{\Delta})\right)
\end{equation*}
and so it is $(m+1)$-fold symmetric, because of \eqref{fi1ms} and
\eqref{fi2ms} with $\Phi_j$ replaced by $\Phi_{j\xi},$ $j=1,2.$

\section{Numerical analysis} In this section we discuss the
numerical analysis of the equation of doubly connected V-states.
There is a number of references on the numerical obtention of
$V$-states (see for instance \cite{DZ} and \cite{DR}).

\subsection{Formulation of the problem}
Recall that a domain $D$ with smooth boundary is a $V$-state if and
only if for some real number $\Omega$, which is the angular velocity
of rotation,
\begin{equation}\label{Vdef1}
\Rea \bigl[ \left(2\Omega\, \overline{z}\, + I(z)\right)\,\vec{\tau}
\bigr] = 0 , \quad z \in
\partial D,
\end{equation}
where  $\vec\tau$ is the unit tangent vector to the boundary of $D$,
positively oriented, and
\begin{equation*}\label{I1}
I(z)= \frac{1}{2 \pi i} \int_{ \partial D} \frac{\overline{\zeta}-
\overline{z}}{\zeta-z}\, d\zeta, \quad z \in \C.
\end{equation*}
If $D$ is doubly connected the boundary has two components, which
are smooth Jordan curves. In the previous sections dealing with
existence issues we have assumed that these curves are of class
$C^{1+\alpha}$ for some $\alpha$ satisfying $0 < \alpha  < 1.$ We
have denoted by $\Gamma_1$ the exterior boundary, and by $\Gamma_2$
the inner boundary. Let us consider proper parameterizations
$z_j(\theta)$, $\theta\in[0, 2\pi],$  of $\Gamma_j, \,j=1,2,$ which
traverse the curves in the counterclockwise direction. Denote by
$z_{j,\theta}$ the derivative of $z_j(\theta)$ with respect to
$\theta.$ Then the single complex equation \eqref{Vdef1} becomes a
system of two real equations

\begin{equation}
\begin{split}
 \Rea[\Big(2 \Omega\, \overline{z_1(\theta)} +
I(z_1(\theta))\Big)z_{1,\theta}(\theta)] = 0,
    \cr
 \Rea[\Big(2 \Omega\,\overline{z_2(\theta)} +
I(z_2(\theta))\Big)z_{2,\theta}(\theta)] = 0.
\end{split}
\end{equation}

Parametrizing the integral defining $I(z)$ this system can be
rewritten as

\begin{align}
\label{e:V-Stateconditions1} \Rea
\bigg[\bigg(2\Omega\overline{z_1(\theta)} & + \frac{1}{2\pi
i}\int_0^{2\pi}\frac{\overline{z_1(\phi)} -
\overline{z_1(\theta)}}{z_1(\phi) -
z_1(\theta)}z_{1,\phi}(\phi)d\phi
    \cr
& - \frac{1}{2\pi i}\int_0^{2\pi}\frac{\overline{z_2(\phi)} -
\overline{z_1(\theta)}}{z_2(\phi) -
z_1(\theta)}z_{2,\phi}(\phi)d\phi\bigg)z_{1,\theta}(\theta)\bigg] =
0,
    \\
\label{e:V-Stateconditions2} \Rea
\bigg[\bigg(2\Omega\overline{z_2(\theta)} & + \frac{1}{2\pi
i}\int_0^{2\pi}\frac{\overline{z_1(\phi)} -
\overline{z_2(\theta)}}{z_1(\phi) -
z_2(\theta)}z_{1,\phi}(\phi)d\phi
    \cr
& - \frac{1}{2\pi i}\int_0^{2\pi}\frac{\overline{z_2(\phi)} -
\overline{z_2(\theta)}}{z_2(\phi) -
z_2(\theta)}z_{2,\phi}(\phi)d\phi\bigg)z_{2,\theta}(\theta)\bigg] =
0.
\end{align}
The second integral in \eqref{e:V-Stateconditions1} and the first
integral in \eqref{e:V-Stateconditions2} are obviously non-singular
(that is, absolutely convergent) because $\Gamma_1$ and $\Gamma_2$
do not intersect. The first integral in \eqref{e:V-Stateconditions1}
and the second integral in \eqref{e:V-Stateconditions2} are also
non-singular, because
\begin{equation}
\label{e:limits} \lim_{\phi \rightarrow
\theta}\frac{\overline{z_j(\phi)} -
\overline{z_j(\theta)}}{z_j(\phi) - z_j(\theta)} =
\frac{\overline{z_{j,\theta}(\theta)}}{z_{j,\theta}(\theta)}, \quad
j=1,2.
\end{equation}
In order to solve the above system it is convenient to work in polar
coordinates
\begin{equation}
\label{e:z1z2} z_1(\theta) = e^{i\theta}\rho_1(\theta), \qquad
z_2(\theta) = e^{i\theta}\rho_2(\theta),
\end{equation}
where $\rho_1$ and $\rho_2$ are given as cosine expansions
\begin{equation}
\label{e:rho1rho2} \rho_1(\theta) = 1 + \sum_{k = 1}^\infty
a_{1,k}\cos(k\,\theta), \qquad \rho_2(\theta) = b + \sum_{k =
1}^\infty a_{2,k}\cos(k\,\theta).
\end{equation}
We are using here that we work the functional space $X$ of section 4
and thus our $V$-states are symmetric with respect to the real axis.
We have normalized so that we get the circle of center the origin
and radius $1$ when all the $a_{1,k}$ vanish and the circle of
center the origin and radius $b$ when all the $a_{2,k}$ vanish. Then
\begin{equation}
\label{e:z1z2cos} z_1(\theta) = e^{i\theta}\left[1 + \sum_{k = 1}^M
a_{1,k}\cos(m\,k\,\theta)\right], \qquad z_2(\theta) =
e^{i\theta}\left[b + \sum_{k = 1}^M
a_{2,k}\cos(m\,k\,\theta)\right],
\end{equation}
and so the problem is reduced to finding numerically the
coefficients $a_{1,k}$ and $a_{2,k}$.  Introducing \eqref{e:z1z2cos}
into \eqref{e:V-Stateconditions1}-\eqref{e:V-Stateconditions2}, we
realize that the errors can be represented as sine expansions of the
form
\begin{equation}
\label{e:V-Stateconditions}
\begin{split}
\Re\bigg[\bigg(2\Omega\overline{z_1(\theta)} & + \frac{1}{2\pi
i}\int_0^{2\pi}\frac{\overline{z_1(\phi)} -
\overline{z_1(\theta)}}{z_1(\phi) -
z_1(\theta)}z_{1,\phi}(\phi)d\phi
    \cr
& - \frac{1}{2\pi i}\int_0^{2\pi}\frac{\overline{z_2(\phi)} -
\overline{z_1(\theta)}}{z_2(\phi) -
z_1(\theta)}z_{2,\phi}(\phi)d\phi\bigg)z_{1,\theta}(\theta)\bigg] =
\sum_{k = 1}^M b_{1,k}\sin(m\,k\,\theta),
    \\
\Re\bigg[\bigg(2\Omega\overline{z_2(\theta)} & + \frac{1}{2\pi
i}\int_0^{2\pi}\frac{\overline{z_1(\phi)} -
\overline{z_2(\theta)}}{z_1(\phi) -
z_2(\theta)}z_{1,\phi}(\phi)d\phi
    \cr
& - \frac{1}{2\pi i}\int_0^{2\pi}\frac{\overline{z_2(\phi)} -
\overline{z_2(\theta)}}{z_2(\phi) -
z_2(\theta)}z_{2,\phi}(\phi)d\phi\bigg)z_{2,\theta}(\theta)\bigg] =
\sum_{k = 1}^M b_{2,k}\sin(m\,k\,\theta),
\end{split}
\end{equation}
where, as before, we take finitely many sines in the error
expansions. Indeed, we choose the same number of cosines and sines.
Therefore, fixed $b$ and $\Omega$, finding a doubly connected
$V$-state is reduced to obtaining a nontrivial root of the nonlinear
equation
\begin{equation}
\mathcal F_{b,\Omega}(a_{1,1}, \ldots, a_{1,M}\,,\, a_{2,1}, \ldots,
a_{2,M}) = (b_{1,1}, \ldots, b_{1,M}\,,\, b_{2,1}, \ldots, b_{2,M});
\end{equation}
where the mapping
$$ \mathcal F_{b,\Omega}\ : \ \mathbb R^{2M}
\longrightarrow \mathbb R^{2M}$$
 is defined from the left hand-side
of \eqref{e:V-Stateconditions} in the obvious way. Notice that we
have trivially $\mathcal F_{b,\Omega}(\mathbf 0) = \mathbf 0$,  for
each value of the parameters $b$ and  $\Omega$. In other words, any
circular annulus is a solution of the problem.


\subsection{Numerical obtention of the $m$-fold $V$-states}

The numerical method that we describe in this section can be applied
with virtually no change to the obtention of simply-connected
$V$-states, and even to more general types of $V$-states.

From the implementation point of view, it is more convenient to work
internally with exponential functions of the form $e^{ik\alpha}$
than with cosines and sines. More precisely, in view of
\eqref{e:z1z2cos} and \eqref{e:V-Stateconditions}, we need the
functions $e^{i(mk+1)\alpha}$, with $k = -M, \ldots, M$. Thus, if we
discretize $[0,2\pi]$ by $N + 1$ equally-spaced nodes $\alpha_j =
2\pi j / N$, $j = 0, \ldots, N$, $N$ has to be chosen for sampling
purposes so that $N \ge 2mM+1$.

All the operations required in \eqref{e:V-Stateconditions}
(obtention of $z_1$ and $z_2$ and their derivatives $z_{1,\alpha}$
and $z_{2,\alpha}$ from the coefficients $a_{1,k}$ and $a_{2,k}$;
and obtention of the coefficients $b_{1,k}$ and $b_{2,k}$) are
computed spectrally via discrete Fourier transforms (DFTs) of $N$
elements, except for the integrals in \eqref{e:V-Stateconditions}
which, bearing in mind \eqref{e:limits}, are numerically evaluated
with spectral accuracy by means of the trapezoidal rule. We choose
$N$ to be a multiple of $m$, $N = m2^r$, so $M = \lfloor
(m2^r-1)/(2m)\rfloor = 2^{r-1}-1$. Then, thanks to the symmetries of
the problem, the DFTs of $N$ elements are reduced to DFTs of $N/m =
2^r$ elements. These DFTs are calculated via the fast Fourier
transform (FFT) algorithm \cite{FJ} in a very efficient way.

In order to find a nontrivial root of $\mathcal F_{b,\Omega}$, we
use a Newton-type iteration. We discretize the
$(2M\times2M)$-dimensional Jacobian matrix $\mathcal J$ of $\mathcal
F_{b,\Omega}$ using just first-order approximations. Fixed $|h|\ll1$
(we have chosen $h = 10^{-9}$), we have
\begin{equation}
\label{e:derivative}
\begin{split}
\frac{\partial}{\partial a_{1,1}} & \mathcal F_{b,\Omega}(a_{1,1},
\ldots, a_{1,M}\,,\, a_{2,1}, \ldots, a_{2,M})
    \cr
& \approx \frac{\mathcal F_{b,\Omega}(a_{1,1} + h, a_{1,2}, \ldots,
a_{1,M}\,,\, a_{2,1}, \ldots, a_{2,M}) - \mathcal
F_{b,\Omega}(a_{1,1}, \ldots, a_{1,M}\,,\, a_{2,1}, \ldots,
a_{2,M})}{h}.
\end{split}
\end{equation}

\noindent Then, the sine expansion of \eqref{e:derivative} gives us
the first row of $\mathcal J$, and so on.

Let us suppose that at the $n$-th iteration we have a good enough
approximation of a root of $\mathcal F_{b,\Omega}$, which we denote
by $(a_{1,1}, \ldots, a_{1,M}\,,\, a_{2,1}, \ldots, a_{2,M})^{(n)}$.
Then, the $(n+1)$-th iteration yields
\begin{equation}
\begin{split}
(a_{1,1}, & \ldots, a_{1,M}\,,\, a_{2,1}, \ldots, a_{2,M})^{(n+1)}
     = (a_{1,1}, \ldots, a_{1,M}\,,\, a_{2,1}, \ldots, a_{2,M})^{(n)}\cr
     &\qquad\qquad -
\mathcal F_{b,\Omega}\left((a_{1,1}, \ldots, a_{1,M}\,,\, a_{2,1},
\ldots, a_{2,M})^{(n)}\right)\cdot [\mathcal J^{(n)}]^{-1},
\end{split}
\end{equation}

\noindent where $[\mathcal J^{(n)}]^{-1}$ denotes the inverse of the
Jacobian matrix corresponding to
$$(a_{1,1}, \ldots, a_{1,M}\,,\,
a_{2,1}, \ldots, a_{2,M})^{(n)}.$$
 This iteration converges in a
small number of steps to a nontrivial root for a large variety of
initial data $(a_{1,1}, \ldots, a_{1,M}\,,\, a_{2,1}, \ldots,
a_{2,M})^{(0)}$. In fact, it is usually enough to perturb the
annulus by assigning a small value to ${a_{1,1}}^{(0)}$ or
${a_{2,1}}^{(0)}$ and leave the other coefficients equal to zero.
Our stopping criterion is
\begin{equation}
\max\left|\sum_{k = 1}^M b_{1,k}\sin(m\,k\,\alpha)\right| <  \rm{tol}
\quad\; \text{and} \quad\; \max\left|\sum_{k = 1}^M
b_{2,k}\sin(m\,k\,\alpha)\right| < \rm{tol},
\end{equation}
where $\rm{tol}= 10^{-12}$, although we get often even smaller errors.

Finally, let us mention that all solutions we obtain by this
procedure satisfy  $a_{1,1}\cdot a_{2,1}<0$. Hence, for coherent
comparisons, we change eventually the sign of all the coefficients
$\{a_{1,k}\}$ and $\{a_{2,k}\}$  in order that, without loss of
generality, $a_{1,1}>0$ and $a_{2,1}<0$.

\begin{figure}[t!]
\center
\includegraphics[width=0.5\textwidth, clip=true]{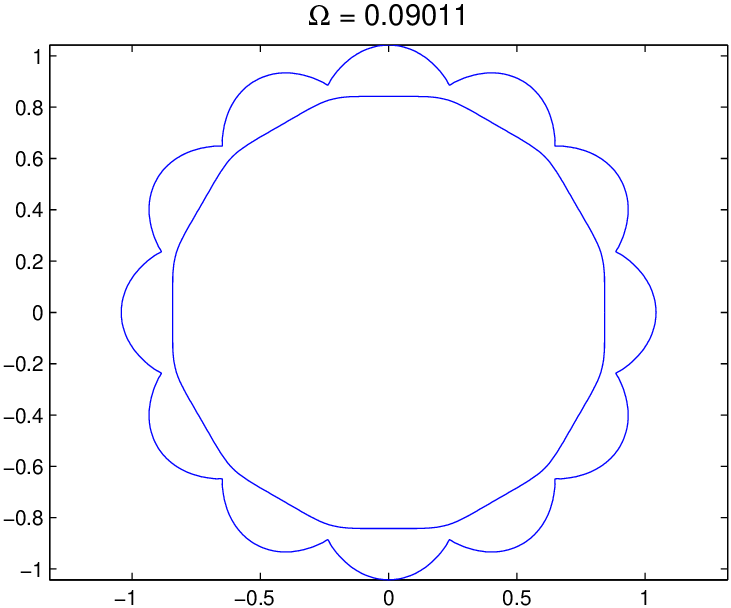}~\includegraphics[width=0.5\textwidth, clip=true]{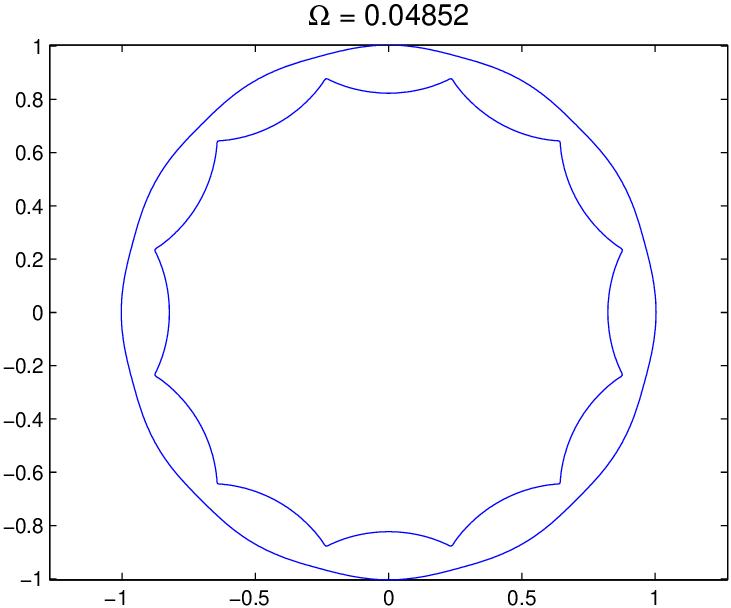}
\caption{Examples of 12-fold $V$-states, for $b = 0.85$.}
\label{f:12-fold}
\end{figure}

In Figure \ref{f:12-fold}, we show two 12-fold $V$-states obtained
via this technique, for $b = 0.85$, using $N = 12\times64=768$
nodes. The left-hand side corresponds to $\Omega = 0.09011$; and the
right-hand side corresponds to $\Omega = 0.04852$. For the
right-hand side, the only initial nonzero coefficient was
${a_{1,1}}^{(0)}=0.06$; and it took nine iterations and about 7.5
seconds to converge. For the left-hand side, the only initial
nonzero coefficient was ${a_{2,1}}^{(0)}=-0.04$; and it took ten
iterations and about 9 seconds to converge. Remark that a couple of
trials may be required until a value of ${a_{1,1}}^{(0)}$ or
${a_{2,1}}^{(0)}$ that enables convergence is found. Once a
$V$-state is found, it can be used as a starting initial value for
finding a new $V$-state with a slightly different $\Omega$ and/or
$b$.

\subsection{Numerical experiments}

According to our main result, Theorem B, given $b\in(0,1)$ the
number of sides $m$ has to be chosen so that
\begin{equation}
f_m(b) = 1 + b^m - \frac{1 - b^2}{2}m < 0.
\end{equation}

\begin{figure}[t!]
\center
\includegraphics[width=0.5\textwidth, clip=true]{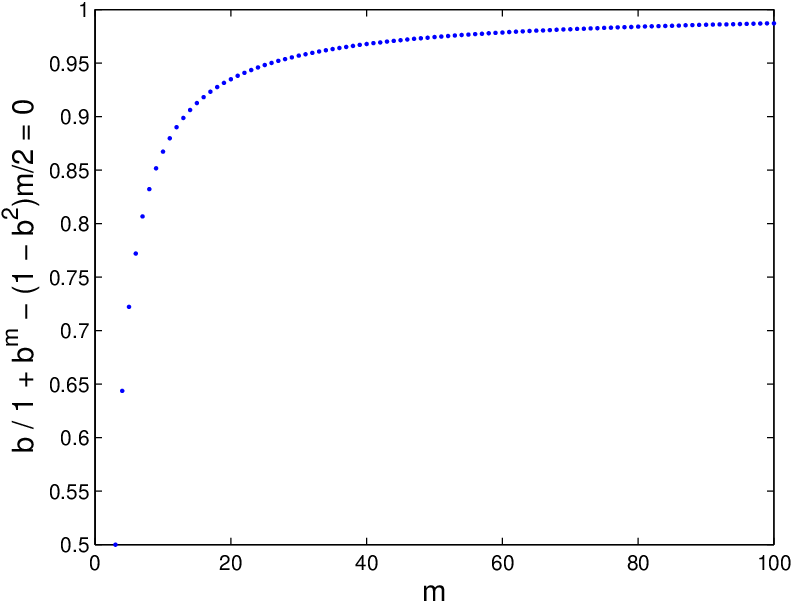}
\caption{Solution of $f_m(b) = 0$, for $m = 3, \ldots, 100$.}
\label{f:rootoff}
\end{figure}

\noindent When $m = 1$ or $m = 2$, $f_m(b)$ is always positive, and
the theorem cannot be applied. When $m \ge 3$, $f_m(0) = 1 - m/2 <
0$ and $f_m(1) = 2 > 0$, so there is at least one $b \in (0,1)$ such
that $f_m(b)=0$. Moreover, since $f_m'(b) = (b + b^{m-1})m
> 0$, $f_{m}$ is a strictly increasing function on the interval
$(0,1)$ and then the equation $f_m(b)=0$ has a unique root on this
interval, which we denote by $b_m$. In Figure \ref{f:rootoff}, we
plot $b_m$ against $m$. The values of $b_m$ have been obtained with
a Newton-type iteration; it is straightforward to check that $b_3 =
1/2$; moreover, $b_m$ tends to $1$ as $m$ grows.

Given $b\in(0, b_m)$, Theorem B guarantees that we can bifurcate
from an annulus with outer radius $1$, inner radius $b$, and angular
velocity $\Omega_m^\pm(b)$, where
\begin{equation}
\label{e:omegamb} \Omega_m^\pm(b) = \frac{1-b^2}{4}\pm
\frac{1}{2m}\sqrt{\left(\frac{m(1-b^2)}{2} - 1\right)^2 - b^{2m}}.
\end{equation}

\noindent Then, on the one hand, $\Omega_m^+(b_m) = \Omega_m^-(b_m)
= (1 - b^2) / 4$;  on the other hand, $\Omega_m^+(0) = (m - 1) /
(2m)$  and $\Omega_m^-(0) = 1 / (2m)$. It is important to remark
that in the analysis of the simply-connected $V$-states of
\cite{DZ}, which corresponds to the limiting case $b = 0$, only
$\Omega_m^+(0)$ appears when bifurcating from a circumference of
radius $1$. This apparently odd behavior will be clarified d in this
section.

In what follows, we take $m = 4$, although everything is immediately
applicable to any $m$. We use always $N = 4\times128 = 512$ nodes.
According to our numerical simulations, there are roughly two
situations: $b$ is ``close'' to $b_4 =\sqrt{\sqrt{2}-1} =
0.6435\ldots$; and $b$ is ``not close'' to $b_4$. We use here
quotation marks because of the informality of the term ``close";
indeed, our aim is to perform a qualitative analysis of 4-fold
$V$-states, rather than a quantitative one.

When $b$ is close enough to $b_4$, it is straightforward to obtain
numerically $V$-states for {each} $\Omega\in(\Omega_m^-,
\Omega_m^+)$ (there is no spectral gap). To illustrate this, we have taken $b = 0.63$.
According to \eqref{e:omegamb}, $\Omega_4^+(0.63) = 0.1674\ldots$
and $\Omega_4^-(0.63) = 0.1341\ldots$; we have calculated the
$V$-states corresponding to the $333$ different values $\Omega =
0.1342, 0.1343, \ldots, 0.1674$  For $\Omega = 0.1342$, the
$V$-state is very close to a circular annulus. Then, as we increase
$\Omega$, the inner boundary resembles more and more a rounded
square; the outer boundary also takes the shape of a rounded square,
rotated of $\pi/4$ degrees with respect to the inner boundary,
although less pronouncedly. However, when $\Omega$ approaches
$\Omega_4^+(0.63)$, we observe the opposite phenomenon, i.e., the
boundaries become more and more circular. For  $\Omega = 0.1674$ we
have again a $V$-state which is very close to a circular annulus.

\begin{figure}[t!]
\center
\includegraphics[width=0.5\textwidth, clip=true]{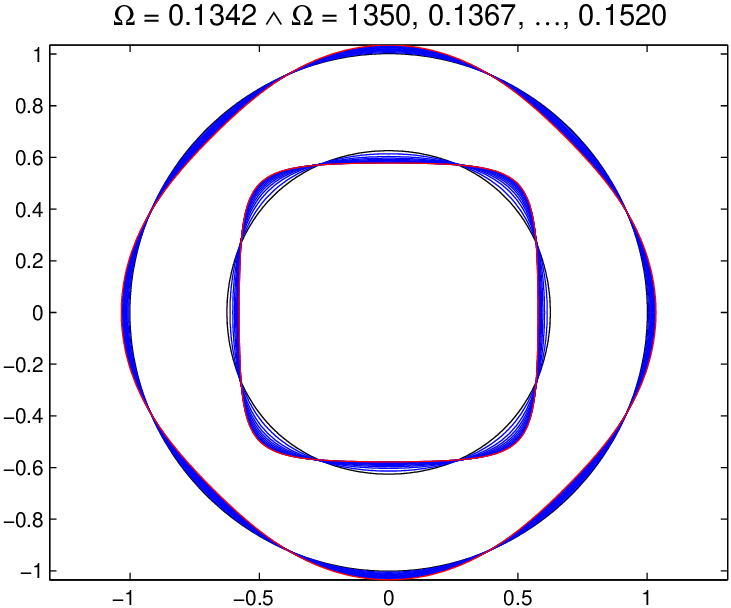}~\includegraphics[width=0.5\textwidth, clip=true]{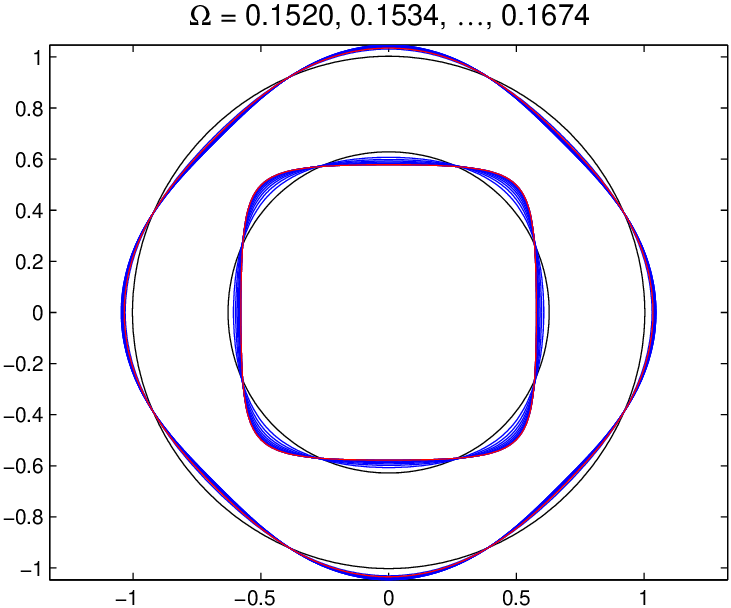}
\caption{Family of 4-fold $V$-states, for $b = 0.63$ and different
$\Omega$.} \label{f:4-fold,b0.63}
\end{figure}
In the left-hand side of Figure \ref{f:4-fold,b0.63}, we have
plotted the $V$-states corresponding to $\Omega = 0.1342$, and to
$\Omega = 0.1350, 0.1367, \ldots, 0.1520$. The $V$-state
corresponding to $\Omega = 0.1342$, in black, is very close to a
circular annulus; while the $V$-state corresponding to $\Omega =
0.1520$, in red, is the $V$-state whose inner boundary is most
pronouncedly a (slightly non-convex) rounded square. In the
right-hand side of Figure \ref{f:4-fold,b0.63}, we have plotted the
$V$-states for $\Omega = 0.1520, 0.1534, \ldots, 0.1674$. The
$V$-state corresponding to $\Omega = 0.1520$ is again in red, while
the $V$-state corresponding to $\Omega = 0.1674$, in black, is very
close to a circular annulus.

\begin{figure}[t!]
\center
\includegraphics[width=0.5\textwidth, clip=true]{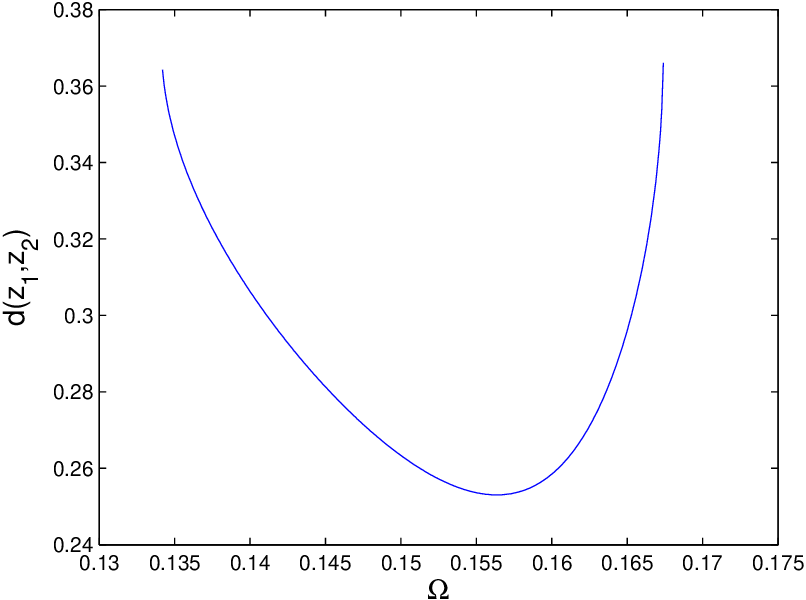}
\caption{Distance between the external and internal boundaries of
the 4-fold $V$ states, for $b = 0.63$ and $\Omega = 0.1342, 0.1343,
\ldots, 0.1674$.} \label{f:distance}
\end{figure}

It is also interesting to compute the distance $d(z_1,z_2)=\inf_{\alpha,\alpha^\prime\in[0,2\pi]}|z_1(\alpha) - z_2(\alpha^\prime)|$
between the boundaries
of  a $V$-state and think of it as a function of  $\Omega$.  This is plotted in
Figure \ref{f:distance}. 
When $\Omega = 0.1342$ and $\Omega = 0.1674$, the distances  respectively $0.3642\ldots$ and $0.3660\ldots$, i.e., they are close
to $1 - b = 0.37$. The minimum distance, $0.2530$, corresponds to
$\Omega = 0.1564$.

However, when $b$ is ``not so close'' to $b_4$, we are able to
obtain 4-fold $V$-states only for $\Omega\in(\Omega_m^-,
\Omega_m^-+\varepsilon^-]$ and
$\Omega\in[\Omega_m^+-\varepsilon^+,\Omega_m^+)$, for certain
$\varepsilon^-$ and $\varepsilon^+$ that depend on $b$. It is
striking that this behavior happens rather soon. Let us take for
instance $b = 0.6$; with $\Omega_4^+(0.6) = 0.1910\ldots$ and
$\Omega_4^-(0.6) = 0.1289\ldots$. When we try to bifurcate from
$\Omega_4^+(0.6)$, we obtain 4-fold $V$-states only until
approximately $\Omega = 0.1755$. In the left-hand side of Figure
\ref{f:4-fold,b0.6}, we have plotted the $V$-states corresponding to
$\Omega = 0.1755$, and to $\Omega = 0.177, 0.179, \ldots, 0.191$.
The $V$-state corresponding to $\Omega = 0.191$, in black, is very
close to a circular annulus. Then, as $\Omega$ gets smaller, the
outer boundary becomes less and less circular, while the inner
boundary resembles more and more to a slightly non-convex rounded
square. At $\Omega = 0,1755$, in red, the inner boundary seems to be
close to developing singularities at the corners of the rounded
square. An analogous situation happens when we try to bifurcate
starting from $\Omega_4^-(0.6)$. We have obtained 4-fold $V$-states
only until approximately $\Omega = 0.158$. In the right-hand side of
Figure \ref{f:4-fold,b0.6}, we have plotted the $V$-states
corresponding to $\Omega = 0.129, 0.132, \ldots, 0.156$, and $\Omega
= 0.158$. The $V$-state corresponding to $\Omega = 0.129$, in black,
is very close to a circular annulus. Then, as $\Omega$ gets larger,
the inner boundary resembles more and more a slightly non-convex
rounded square, while the outer boundary, unlike in the previous
case, remains always rather close to a circumference. At $\Omega =
0.158$, in red, the inner boundary seems to be close to developing
singularities at the corners of the rounded square.
\begin{figure}[t!]
\center
\includegraphics[width=0.5\textwidth, clip=true]{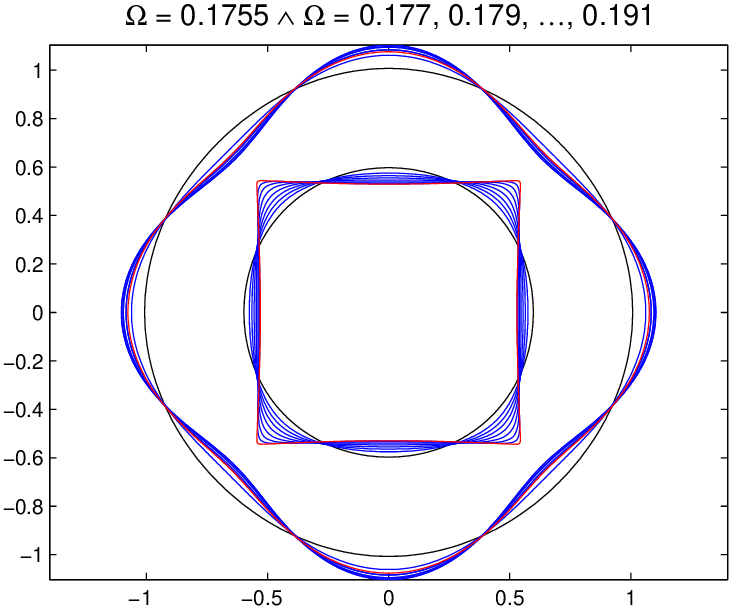}~\includegraphics[width=0.5\textwidth, clip=true]{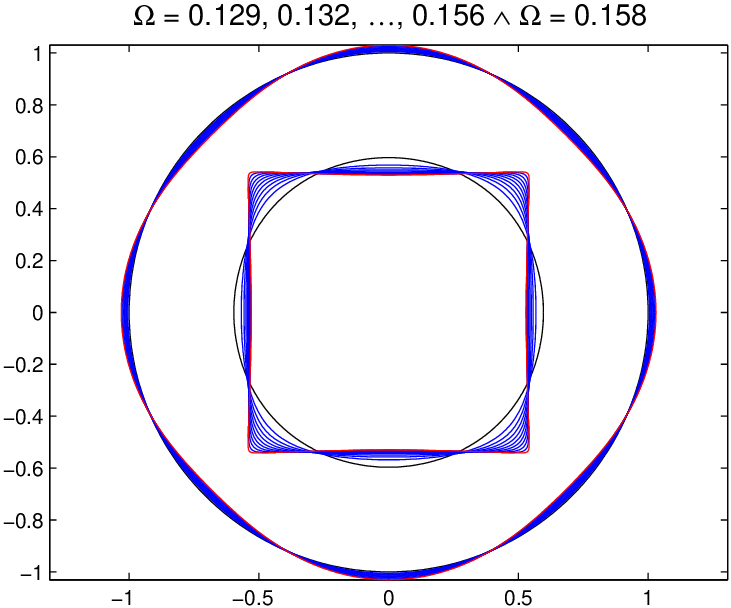}
\caption{Family of 4-fold $V$-states, for $b = 0.6$ and different
$\Omega$. In the left-hand side, we have started to bifurcate from
$\Omega_4^+(0.6)=0.1910\ldots$; while, in the right-hand side, we
have started to bifurcate from $\Omega_4^-(0.6) = 0.1289\ldots$.}
\label{f:4-fold,b0.6}
\end{figure}

Summarizing, for $\Omega\in[\Omega^-(0.6)+\varepsilon^-,
\Omega^+(0.6)-\varepsilon^+]$, where
$\Omega^-(0.6)+\varepsilon^-\approx 0.158$ and
$\Omega^+(0.6)-\varepsilon^+\approx0.1755$, numerical instabilities
appear and we are unable to obtain bifurcated $V$-states. Remark
that something similar happens with the examples of the 12-fold
$V$-states in Figure \ref{f:12-fold}, which are also limiting cases;
in fact, the singularities are even more evident in that figure. It
is also worth mentioning that the boundaries of the left-hand side
of Figure \ref{f:12-fold} are very close from each other at some
points. Furthermore, by choosing carefully the parameters, it is
possible to find $V$-states whose boundaries seem almost to touch
each other.

We have also computed $V$-states for smaller $b$. In Figure
\ref{f:4-fold,b0.4}, we have taken $b = 0.4$; with $\Omega_4^+(0.4)
= 0.2949\ldots$ and $\Omega_4^-(0.4) = 0.1250\ldots$. When we start
to bifurcate from $\Omega_4^+(0.4)$, the inner boundaries almost do
not change and remain close to a circumference all the time, while
the outer boundaries get closer and closer to a non-convex rounded
square. In the left-hand side of Figure \ref{f:4-fold,b0.4}, we have
plotted the $V$-states corresponding to $\Omega = 0.267, 0.270,
\ldots, 0.294$, and to $\Omega = 0.2949$. The $V$-state
corresponding to $\Omega = 0.2949$, in black, is very close to a
circular annulus, while the $V$-state corresponding to $\Omega =
0.267$, in red, seems to be close to developing singularities. We
have exactly the opposite situation when we start to bifurcate from
$\Omega_4^-(0.4)$, because the outer boundaries are the ones that
remain close to a circumference, while the inner boundaries tend to
a slightly non-convex rounded square. In the right-hand side of
Figure \ref{f:4-fold,b0.4}, we have plotted the $V$-states
corresponding to $\Omega = 0.126, 0.130, \ldots, 0.146$. The
$V$-state for $\Omega = 0.146$, in black, is very close to a
circular annulus, while the $V$-state corresponding to $\Omega =
0.126$, in red, seems to be close to developing singularities.

\begin{figure}[t!]
\center
\includegraphics[width=0.5\textwidth, clip=true]{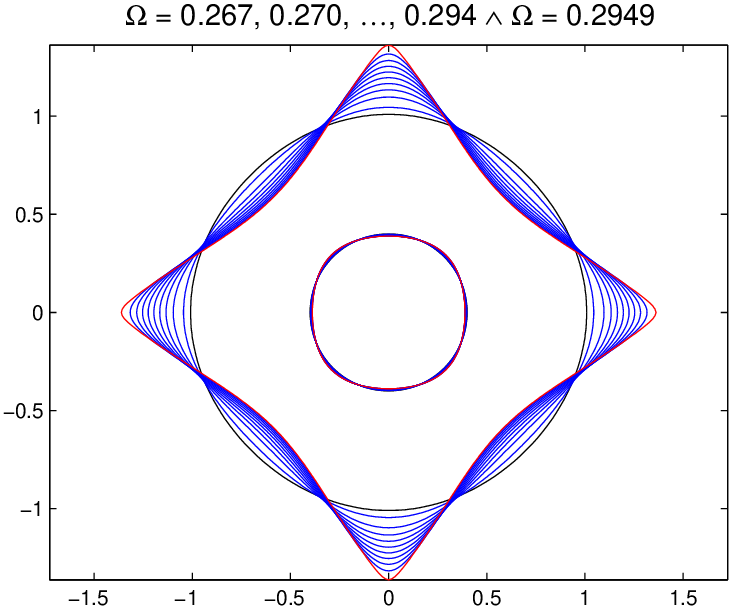}~\includegraphics[width=0.5\textwidth, clip=true]{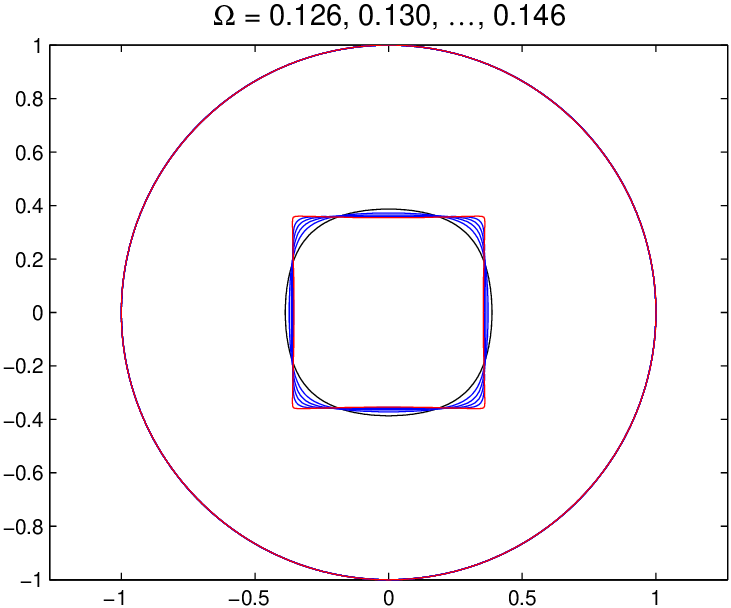}
\caption{Family of 4-fold $V$-states, for $b = 0.4$ and different
$\Omega$. In the left-hand side, we have started to bifurcate from
$\Omega_4^+(0.4)=0.2949\ldots$; while, in the right-hand side, we
have started to bifurcate from $\Omega_4^-(0.4) = 0.1250\ldots$.}
\label{f:4-fold,b0.4}
\end{figure}

All the conclusions for $b = 0.4$ are valid for smaller $b$,
although even more exaggerated, as is clear from Figure
\ref{f:4-fold,b0.2}, where we have taken $b = 0.2$; with
$\Omega_4^+(0.2) = 0.3549\ldots$ and $\Omega_4^-(0.2) =
0.1250\ldots$. In fact, all the previous considerations apply, so we
do not mention them again. Furthermore, Figures \ref{f:4-fold,b0.4}
and \ref{f:4-fold,b0.2} explain the apparently odd behavior
mentioned above, when we pointed that in the doubly connected case
we could bifurcate from the annulus at two values of $\Omega$, while
in the simply-connected case, there was only one such value. Indeed,
when $b$ tends to 0, the $V$-states obtained after bifurcating from
$\Omega_4^-(b)$ just tend to the unit circle, while those obtained
after bifurcating from $\Omega_4^+(b)$ tend to a simply connected
4-fold $V$-state.

\begin{figure}[t!]
\center
\includegraphics[width=0.5\textwidth, clip=true]{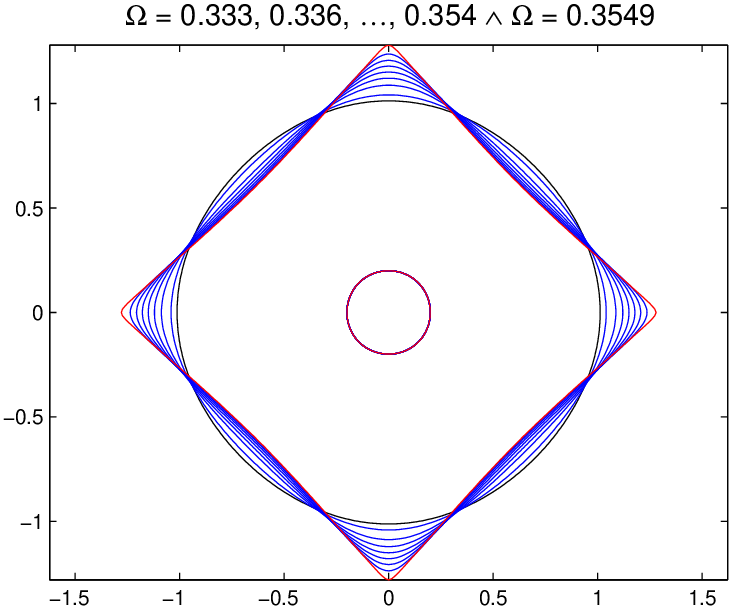}~\includegraphics[width=0.5\textwidth, clip=true]{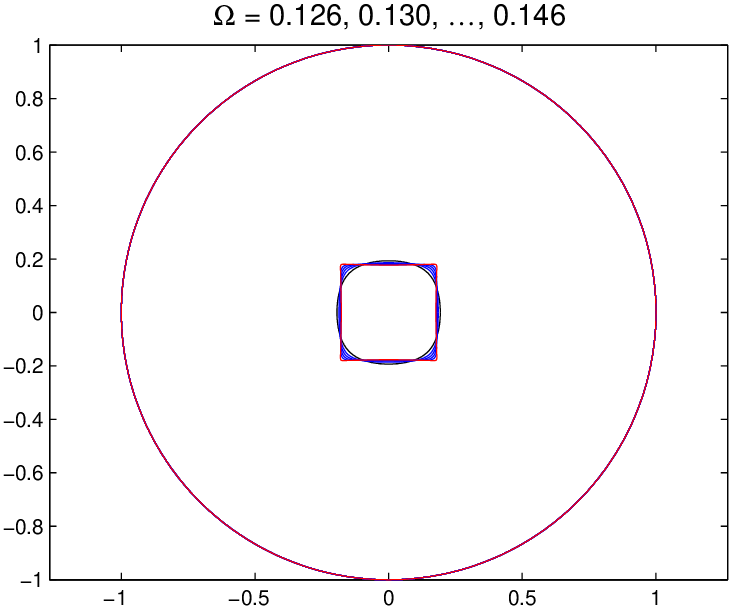}
\caption{Family of 4-fold $V$-states, for $b = 0.2$ and different
$\Omega$. In the left-hand side, we have started to bifurcate from
$\Omega_4^+(0.2)=0.3549\ldots$; while, in the right-hand side, we
have started to bifurcate from $\Omega_4^-(0.2) = 0.1250\ldots$.}
\label{f:4-fold,b0.2}
\end{figure}


\section{Conclusion} We have shown that simple eigenvalues $\lambda$
are obtained by requiring that $\Delta_m(\lambda,b)=0$  (see
\eqref{DSn}) for some frequency $ 0 \le m \neq 1$ or, for $m=1,$  by
$\lambda= (1+b^2)/2$ with $b$ not belonging to the sequence $\{ b
\in (0,1): \varphi_n(b)=0, \; \text{for some} \; n \ge 2 \}$, where
$\varphi_n$ is defined in \eqref{E3}. If $b$ belongs to this
sequence, then $\lambda= (1+b^2)/2$ is a double eigenvalue. One can
solve the equation $\Delta_m(\lambda,b)=0$ for $\lambda$ and then
compute the angular velocity of rotation $\Omega=(1-\lambda)/2$ .
One gets the formula
\begin{equation*}\label{angular}
\Omega_m = \frac{1-b^2}{4}\pm
\frac{1}{2(m+1)}\sqrt{\bigl(\frac{(m+1)(1-b^2)}{2}-1\bigr)^2
-b^{2(m+1)}},
\end{equation*}
which should be compared to \cite[(4.1), p. 162]{DR}. The eigenvalue
$\sigma(m,a)$ found in \cite{DR} is exactly $m\, \Omega_{m-1}$ with $b$
replaced by $a.$

We have proved that there exists a curve of non-annular $V$-states
that bifurcates from the annulus $\{z: b <|z| <1\}$ for all
eigenvalues associated with frequencies $m \ge 2.$ Given $b \in
(0,1),$ if the frequency $m$ satisfies $m \ge (3+b^2)/(1-b^2)$ then,
by \eqref{gen3}, the equation $\Delta_m(\lambda,b)=0$ has two real
solutions which are simple eigenvalues at which one can bifurcate.
Thus, given any annulus of the form $A=\{z: b <|z| <1\}$, there are
non-annular $V$-states bifurcating at $A.$ They are $(m+1)$-fold
symmetric as we proved in subsection 8.3. This adds a valuable
detailed information to the concise statement of Theorem A and
proves the more precise statement of Theorem B.

There are two simple eigenvalues for which all available criteria
for bifurcation we have found in the literature fail.  These are
$\lambda=b^2$ and $\lambda= (1+b^2)/2$ with $b$ not belonging to the
sequence $\{ b \in (0,1): \varphi_n(b)=0, \; \text{for some} \; n
\ge 2 \}$. More precisely, the transversality condition $(d)$ in Crandall-Rabinowitz's Theorem is not satisfied. For these eigenvalues we do not have any argument ad hoc
to show that bifurcation is possible, nor we have an argument to
show that bifurcation cannot happen. Deciding whether bifurcation
takes place at these simple eigenvalues remains an open question.

\begin{gracies}
We thank L.Vega for posing the problem, P. Luzzatto-Fegiz for
bringing to our attention the paper \cite{DR} and D.R. Dritschel for some very useful
correspondence.  This work was
partially supported by the grants  2014SGR75 (Generalitat de
Catalunya), IT641-13 (Basque Goverment), MTM2011-24054 and
MTM2013-44699 (Ministerio de Econom\'{\i}a y Competividad)  and the ANR
project Dyficolti ANR-13-BS01-0003-01.
\end{gracies}

\vspace*{.55cm}

\begin{tabular}{l}
Taoufik Hmidi\\
IRMAR, Universit\'e de Rennes 1\\
Campus de Beaulieu, 35042 Rennes cedex, France \\
{\it E-mail:} {\tt thmidi@univ-rennes1.fr}\\ \\
Francisco de la Hoz\\
Department of Applied Mathematics, Statistics and Operations
Research \\
University of the Basque Country UPV-EHU\\
 48940 Leioa, Spain \\
{\it E-mail:} {\tt francisco.delahoz@ehu.es}
 \\ \\
Joan Mateu and Joan Verdera\\
Departament de Matem\`{a}tiques\\
Universitat Aut\`{o}noma de Barcelona\\
08193 Bellaterra, Barcelona, Catalonia\\
{\it E-mail:} {\tt mateu@mat.uab.cat}\\
{\it E-mail:} {\tt jvm@mat.uab.cat}
\end{tabular}
\end{document}